\documentclass[10pt]{article}%

\usepackage[utf8]{inputenc}

\usepackage{amsmath}

\usepackage{amsfonts}
\usepackage{mathrsfs}
\usepackage{amssymb, color}
\usepackage[linkcolor=black,anchorcolor=black,citecolor=black]{hyperref}
\usepackage{graphicx}
\numberwithin{equation}{section}
\usepackage[body={15.5cm,21cm}, top=3cm]{geometry}%
\providecommand{\U}[1]{\protect\rule{.1in}{.1in}}
\providecommand{\U}[1]{\protect \rule{.1in}{.1in}}
\newtheorem{theorem}{Theorem}[section]

\newtheorem{definition}[theorem]{Definition}

\newtheorem{lemma}[theorem]{Lemma}

\newtheorem{proposition}[theorem]{Proposition}
\newtheorem{remark}[theorem]{Remark}

\newtheorem{assumption}[theorem]{Assumption}
\newenvironment{proof}[1][Proof]{\noindent \textbf{#1.} }{\  \rule{0.5em}{0.5em}}
\DeclareMathOperator*{\esssup}{ess\,sup}

\def \E{\mathsf{E}}

\def \P{\mathsf{P}}

\def \E{\mathsf{E}}
\def \P{\mathsf{P}}

\begin{document}
	\title{Conditional Reflected Backward Stochastic Differential Equations}
	\author{ 	Hanwu Li\thanks{Research Center for Mathematics and Interdisciplinary Sciences, Shandong University, Qingdao 266237, Shandong, China. lihanwu@sdu.edu.cn.}
	\thanks{Frontiers Science Center for Nonlinear Expectations (Ministry of Education), Shandong University, Qingdao 266237, Shandong, China.}
    \thanks{Shandong Province Key Laboratory of Financial Risk, Shandong University, Qingdao 266237, Shandong, China.}}
	\date{}
	\maketitle
	
	\begin{abstract}
     In this paper, we introduce a new type of reflected backward stochastic differential equations (reflected BSDEs), referred to as conditional reflected BSDEs, in which the conditional expectation of the $Y$-component of the solution with respect to a given sub-$\sigma$-algebra appears in both the generator and the obstacle. Under a Lipschitz condition on the loss function with  sufficiently small Lipschitz constants, we establish the existence and uniqueness of the solution for both the case of Lipschitz generator and the case of quadratic generator, bounded terminal value and bounded loss function. When the generator has quadratic growth and the terminal value and loss function are unbounded, we establish the uniqueness with the help of $\theta$-method. Finally, a penalization approximation scheme holds under an additional monotonicity structure imposed on the generator and the loss function, without any restriction on the Lipschitz constants of the loss function.  
	\end{abstract}


    \textbf{MSC-classification}: 60H10

\section{Introduction}   

  Backward stochastic differential equation (BSDE for short), whose solution is a pair of adapted processes $(Y,Z)$ evolving according to the following equaiton
\begin{align}\label{BSDE}
    Y_t=\xi+\int_t^T f(s,Y_s,Z_s)ds-\int_t^T Z_sdB_s,
\end{align}
has been extensively studied since the pioneering work by Pardoux and Peng \cite{PP}.
It has been widely recognized that the theory of BSDEs serves as a powerful tool for investigating numerous problems in stochastic control, stochastic representations of partial differential equations (PDEs) and engineering, investment science including mathematical finance, differential games and insurance (see, e.g., \cite{DE,EPQ,H98,Norberg,PP92}). Since its inception, extensive extensions have been developed in various directions. One branch is the so-called mean-field BSDE,  where the distribution of its solution also enters into the generator, which is motivated by the models of large stochastic particle systems with mean-field interactions, as well as mean-field limits in diverse fields such as statistical mechanics and physics, quantum mechanics, and quantum chemistry (see, e.g., \cite{BDLP,BLP,HHTW}). On the other hand, the need to impose additional conditions on stochastic control problems and mathematical finance has naturally prompted the investigation of constrained BSDEs. More precisely, El Karoui et al. \cite{KKPPQ} first tackled the case where the solution $Y$ is required to remain above a prescribed obstacle $L$, that is, 
\begin{align}\label{constraint}
    Y_t\geq L_t, \ t\in[0,T].
\end{align}
In order to satisfy this constraint, a non-decreasing process $K$, regarded as the force aiming to push the solution upward, is added into equation \eqref{BSDE}, i.e.,
\begin{align}\label{RBSDE}
    Y_t=\xi+\int_t^T f(s,Y_s,Z_s)ds-\int_t^T Z_sdB_s+K_T-K_t.
\end{align}
By the comparsion theorem for BSDEs, the larger the compensating term $K$, the larger the associated  process $Y$. It is natural to find the minimal solution $Y$ satisfying constraint \eqref{constraint}, which can be guaranteed by the following Skorohod condition
\begin{align}\label{Skorohod}
    \int_0^T(Y_t-L_t)dK_t=0.
\end{align}
The class of constrained BSDEs has been widely generalized along various directions, including the doubly reflected case \cite{CK,HH}, the multi-dimensional case \cite{GP96,HZ,HT}, the case of mean reflection \cite{BEH,FS,Li} and numerous applications in  investment strategies under portfolio constraints \cite{CKS,PX}, switching problems \cite{HJ}, optimal stopping \cite{CR} and many others. We refer to the recent paper \cite{PR} for further discussions and references.

Recently, Djehiche et al. \cite{DEH} investigated the reflected BSDEs of mean-field type taking the following form
\begin{equation*}
\begin{cases}
Y_t=\xi+\int_t^T f(s,Y_s,\P_{Y_s},Z_s)ds-\int_t^T Z_s dB_s+K_T-K_t, \\
Y_t\geq h(t,Y_t,\P_{Y_t}), \ t\in[0,T],\\
\int_0^T \left(Y_t- h(t,Y_t,\P_{Y_t})\right)dK_t=0,
\end{cases}
\end{equation*}
where $\P_{Y_t}$ represents the probability distribution of the process $Y$ at time $t$. Under Lipschitz and intergrability conditions on the generator $f$ and the loss function $h$ with the Lipschitz constants of $h$ being sufficiently small, they establish the existence and uniqueness result via the Snell envelope method and fixed-point argument when the generator $f$ does not depend on the $z$-term. For the general case, the authors use the penalization approach to prove existence and uniqueness of the solution, which requires additional monotonicity assumptions on the generator and obstacle with respect to the measure component.  Then, Djehiche et al. \cite{DDZ} and Hu et al. \cite{HMW} apply a new method with the help of $g$-expectations and the nonlinear optimal stopping problems, which allows to get rid of the monotonicity assumption in \cite{DEH}. They extend the well-posedness of the mean-field reflected BSDE to the setting when the driven process has jumps and to the case where the generator has quadratic growth, respectively. The readers may refer to \cite{CHM,DD} for the mean-field doubly reflected BSDEs.

Motivated by financial problems involving state constraints under partial information, such as American option pricing and recursive reflected utility maximization with partial information, Hu, Huang and Li \cite{HHL} proposed the following BSDE with conditional reflection
\begin{equation*}
\begin{cases}
Y_t=\xi+\int_t^T f(s,Y_s,Z_s)ds-\int_t^T Z_s dB_s+K_T-K_t, \\
\E[Y_t|\mathcal{G}_t]\geq \E[S_t|\mathcal{G}_t], \ t\in[0,T], \\
K\in \mathcal{A}^2_{\mathbb{G}} \textrm{ and }
\int_0^T \E[Y_t-S_t|\mathcal{G}_t]dK_t=0.
\end{cases}
\end{equation*}
Here, $\mathbb{G}=\{\mathcal{G}_t\}_{t\in[0,T]}$ is a given  subfiltration and $\mathcal{A}^2_{\mathbb{G}}$ is the collection of $\mathbb{G}$-adapted, nondecreasing and continuous function starting from the origin. The authors obtain its existence and uniqueness under mild conditions on the parameters $\xi,f,S$ and usual conditions as well as left-quasi-continuity of $\mathbb{G}$ by combining the Snell envelope method with the Skorohod lemma. Subsequently, Li \cite{Li26} considered the BSDE with two conditional reflections and provided an application to the starting and stopping problem in reversible investments under partial information.


In this paper, the objective is to investigate the following kind of conditional reflected BSDE 
\begin{equation}\label{CRBSDE}
\begin{cases}
Y_t=\xi+\int_t^T f(s,Y_s,\E[Y_s|\mathcal{G}_s],Z_s)ds-\int_t^T Z_s dB_s+K_T-K_t, \\
Y_t\geq h(t,Y_t,\E[Y_t|\mathcal{G}_t]), \ t\in[0,T],\\
\int_0^T \left(Y_t- h(t,Y_t,\E[Y_t|\mathcal{G}_t])\right)dK_t=0,
\end{cases}
\end{equation}
where both the generator $f$ and the loss function $h$  depend on the conditional expectation of the first component of its solution. When the subfiltration is trivial, the conditional reflected BSDE degenerates into the mean-field reflected BSDE investigated in \cite{DEH}. As shown in \cite{DD,DEH}, one of the applications of the mean-field reflected BSDE is pricing life insurance contracts  with surrender options. Let $(Y^{\ell,N},Z^{\ell,N})$ be the characteristics of the prospective reserve of each homogeneous life insurance policy, $\ell=1,2,\cdots,N$. An implicit hypothesis in \cite{DD,DEH} is that for each $t\in[0,T]$, the sequence $\{Y^{\ell,N}_t\}_{\ell=1}^N$ is independent and identically distributed. If it is a sequence of conditionally independent \footnote{As claimed in \cite{Smyth}, conditional independence constitutes a highly useful construct within practical modeling frameworks. While such independence assumptions may not strictly hold in real-world scenarios, they can be reasonably imposed to construct approximate models. This treatment is conceptually analogous to the practice of approximating nonlinear functions via linear first-order approximations in applied mathematical modeling. 
Typical examples of conditional independent random variables can also be found in \cite{Smyth}. Besides, the notion of conditional independence shares key properties of the notion of graph separation (see \cite{GP,SB} for more details).} random variables with the same conditional distribution (w.r.t. a sub-$\sigma$-field $\mathcal{G}_t$), we will encounter the conditional expectation appearing in both the generator and the solvency level  of the nonlinear reserving model (more details can be found in Remark \ref{example} below). That is the reason why we investigate conditional reflected BSDEs. In order to establish the existence and uniqueness result, motivated by \cite{DDZ,HMW}, the key point is the following nonlinear Snell envelope representation for the reflected BSDE \eqref{constraint}-\eqref{Skorohod}:
\begin{align*}
    Y_t=\esssup_{\tau \textrm{ stopping time } \geq t}\mathcal{E}^g_{t,\tau}[\xi I_{\{\tau=T\}}+L_\tau I_{\{\tau<T\}}],
\end{align*}
where $\mathcal{E}^{g}_{t,\tau}[\cdot]$ is the $g$-expectation induced by the solution to a certain BSDE. Then, by constructing a contraction mapping, we obtain the well-posedness of conditional reflected BSDEs for both the case of Lipschitz generator and the case of quadratic generator with bounded terminal value and bounded loss function. For the case where the generator has quadratic growth and the terminal value and loss function are unbounded, with the help of $\theta$-method proposed in \cite{BH2}, we derive the uniqueness of the solution. Besides, the comparison theorem for the conditional reflected BSDEs \eqref{CRBSDE} is firstly established under some additional monotonicity assumptions for both the generator and the loss function. To the best of our knowledge, the result of this comparison theorem is entirely new, even for mean-field reflected BSDEs (with $\mathcal{G}_t=\mathcal{F}_0$, $t\in[0,T]$ in \eqref{CRBSDE}).

It is worth pointing out that, when investigating the mean-field reflected BSDEs and the conditional reflected BSDEs by the contraction mapping method, there is an additional requirement that the Lipschitz constants of the loss function are sufficiently small (see Eqs. (18), (28) in \cite{CHM}, Eq. (2.6) in \cite{DD}, Eq. (2.4) in \cite{DDZ}, Eq. (3.8) in \cite{DEH}, Eqs. (9), (14) and (17) in \cite{HMW} and Eqs. \eqref{gamma1gamma2'} and \eqref{gamma1gamma2''}). Another contribution of the present paper is that we may get rid of these hypotheses by combining the penalization method and the monotone convergence for $g$-supersolutions (see \cite{P99}). Actually, a penalization-type constructive argument has been used in \cite{DEH}, which may also reinterpret as a recursive scheme where, at step $n$, the solution to the BSDE $Y^n$ is reflected on the obstacle generated by the solution $Y^{n-1}$. Compared with their construction, we use a purely penalization approach and we do not need to assume that the loss function is non-decreasing in $y$ as \cite{DEH}. However, the penalization method comes at the cost that the solution can no longer be guaranteed to be continuous while satisfying the Skorohod condition. Fortunately, the Skorohod condition generally ensures the minimality of the solution, which is also preserved for solutions constructed by approximation via penalization.  %

The paper is organized as follows. In Section 2, we recall some notions and results on BSDEs and conditional expectation BSDEs. Then, we investigate the conditional reflected BSDEs with Lipschitz generator in Section 3 and establish the existence and uniqueness as well as the comparison theorem. The case when the generator has quadratic growth is studied in Section 4. Finally, utilizing the penalization method, we construct the solutions to the conditional reflected BSDEs in Section 5.

\section{Preliminaries}

Consider a complete probability space $(\Omega,\mathcal{F},\P)$  under which $B$ is a standard one-dimensional Brownian motion (actually, all the results in this paper still hold for the case driven by $d$-dimensional Brownian motion).  Let $\mathbb{F}=\{\mathcal{F}_t\}_{t\in[0,T]}$ be the natural filtration generated by $B$ augmented by all $\P$-null set. Let $\E[\cdot]$ (resp., $\E[\cdot|\mathcal{F}_t]$) be the expectation (resp., conditional expectation w.r.t. $\mathcal{F}_t$) under $\P$. For each $p\geq 1$, given $\beta>0$ and $0\leq s\leq t\leq T$, we first introduce the following spaces which are frequently used in this paper:

\begin{itemize}
    \item $L^p(\mathcal{F}_t)$: the set of real-valued $\mathcal{F}_t$-measurable random variable $\xi$ such that $\E[|\xi|^p]<\infty$.
    \item $L^\infty(\mathcal{F}_t)$: the set of real-valued $\mathcal{F}_t$-measurable random variable $\xi$ such that $\esssup_{\omega\in\Omega}|\xi(\omega)|<\infty$.
    \item $\mathbb{L}^p(\mathcal{F}_t)$: the set of real-valued $\mathcal{F}_t$-measurable random variable $\xi$ such that $\E[\exp(p|\xi|)]<\infty$.
    \item $\mathcal{T}_{s,t}$:  the set of all $\mathbb{F}$-stopping times $\tau$ such that $\tau\in[s,t]$, a.s.
    \item $\mathcal{S}^p_{\beta,c}(s,t)$: the set of real-valued, continuous and adpated processes $Y$, such that $$\|Y\|^p_{\mathcal{S}^p_\beta(s,t)}:=\E\left[\sup_{r\in[s,t]}e^{\beta r}|Y_r|^p\right]<\infty.$$
    \item $\mathcal{S}^\infty_c(s,t)$: the set of real-valued, continuous and adpated processes $Y$, such that $$\|Y\|_{\mathcal{S}^\infty(s,t)}:=\esssup_{(t,\omega)\in [s,t]\times\Omega}|Y_t(\omega)|<\infty.$$
    \item $\mathbb{S}^p_c(s,t)$: the set of processes $Y$ such that $e^Y\in \mathcal{S}^p_{0,c}(s,t)$.
    \item $\mathcal{S}^p_{\beta,ci}(s,t)$: the subset of $\mathcal{S}^p_{\beta,c}(s,t)$ such that the process $K$ is non-decreasing and $K_s=0$.
    \item $\mathcal{S}_{ci}(s,t)$: the set of all processes $K\in \mathcal{S}^p_{0,ci}(s,t)$ for any $p\geq 1$.
    \item $\mathcal{H}^p(s,t)$: the set of real-valued, progressively measurable processes $Z$ such that 
    \begin{align*}
        \|Z\|^p_{\mathcal{H}^p(s,t)}:=\E\left[\left(\int_s^t |Z_r|^2dr\right)^{\frac{p}{2}}\right]<\infty.
    \end{align*}
    \item $\mathbb{L}(\mathcal{F}_t)$: the collection of all random variables $\xi\in \mathbb{L}^p(\mathcal{F}_t)$ for any $p\geq 1$, and $\mathcal{H}, \mathbb{S}_c$ are defined similarly.
    \item $BMO(s,t)$:  the set of real-valued, progressively measurable processes $Z$ such that 
    \begin{align*}
        \|Z\|^2_{BMO(s,t)}:=\sup_{\tau\in \mathcal{T}_{s,t}}\esssup_{\omega\in\Omega}\E\left[\int_\tau^t |Z_r|^2dr\Big|\mathcal{F}_\tau\right]<\infty.
    \end{align*}
    \item $\mathcal{S}^p(s,t)$: the set of real-valued, c\`{a}dl\`{a}g and adpated processes $Y$ such that  $\E\left[\sup_{r\in[s,t]}|Y_r|^p\right]<\infty$. 
    \item $\mathcal{S}^p_{i}(s,t)$: the subset of $\mathcal{S}^p(s,t)$ such that the process $K$ is non-decreasing and $K_s=0$.
\end{itemize}
   We set $\mathcal{S}^p_c(s,t)=\mathcal{S}^p_{0,c}(s,t)$, $\mathcal{S}^p_{ci}(s,t)=\mathcal{S}^p_{0,ci}(s,t)$. If $s=0,t=T$, we omit the brankets. Clearly, we have $\mathcal{S}^2_{\beta,c}=\mathcal{S}^2_c$ and $\mathcal{S}^2_{\beta,ci}=\mathcal{S}^2_{ci}$.

We first recall some basic notation and results on BSDE and $g$-expectation, where $g$ is the generator of a BDSE satisfying the following usual conditions.

\begin{assumption}\label{assg}
    The generator $g$ is a map from $\Omega\times[0,T]\times \mathbb{R}\times\mathbb{R}$ to $\mathbb{R}$ such that 
\begin{itemize}
    \item [(a)] for any $y,z\in\mathbb{R}$, $\{g(t,y,z)\}_{t\in[0,T]}$ is progressively measurable and $\{g(t,0,0)\}_{t\in[0,T]}\in\mathcal{H}^2$,
    \item [(b)] there exists a constant $L>0$ such that $\P$-a.s. for any $t\in[0,T]$ and any $y_i,z_i\in\mathbb{R}$, $i=1,2$
\begin{align*}
|g(t,y_1,z_1)-g(t,y_2,z_2)|\leq L(|y_1-y_2|+|z_1-z_2|).
\end{align*}
\end{itemize}
\end{assumption}

\begin{definition}
    Given $\tau\in \mathcal{T}_{0,T}$, $\xi\in L^2(\mathcal{F}_\tau)$ and $g$ satisfying Assumption \ref{assg}, let $(Y^{\tau,\xi},Z^{\tau,\xi})$ be the solution to BSDE with terminal time $\tau$, terminal value $\xi$ and generator $g$. For any stopping time $\sigma\leq \tau$, the conditional $g$-expectation $\mathcal{E}^g_{\sigma,\tau}:L^2(\mathcal{F}_\tau)\rightarrow L^2(\mathcal{F}_\sigma)$ is defined as follows:
    \begin{align*}
        \mathcal{E}^g_{\sigma,\tau}[\xi]=Y^{\tau,\xi}_\sigma.
    \end{align*}
\end{definition}

Motivated by the a priori estimates for BSDEs driven by $G$-Brownian motion (see the proof of Proposition 3.9 in \cite{HJPS}), We have the following a priori estimates for BSDEs. 

\begin{proposition}\label{a priori estimates}
   Given $\tau\in \mathcal{T}_{0,T}$ and $\xi^i\in L^2(\mathcal{F}_\tau)$, $i=1,2$, let $f^1$ be a generator with Lipschitz constant $C$ and $f_2$ be a generator. Let  $(Y^i,Z^i)$ be the solution to BSDE associated with terminal time $\tau$, terminal value $\xi^i$ and generator $f^i$, $i=1,2$. For any $t\in[0,\tau]$, set $\delta{Y}_t=Y^1_t-Y^2_t$, $\delta{\xi}=\xi^1-\xi^2$ and 
    \begin{align*}
        \delta_2{f}_t=f^1(t,Y^2_t,Z^2_t)-f^2(t,Y^2_t,Z^2_t).
    \end{align*}
    For any $\alpha,\beta$ such that $\alpha\in(1,2]$ and $\beta>\alpha C+\alpha-1+\frac{\alpha C^2}{\alpha-1}$, we have  
    \begin{align*}
        e^{\beta t}|\delta Y_t|^\alpha\leq \E\left[e^{\beta \tau}|\delta \xi|^\alpha|\mathcal{F}_t\right]+\E\left[\int_t^\tau e^{\beta s}|\delta_2 f_s|^\alpha ds\Big|\mathcal{F}_t\right].
    \end{align*}
\end{proposition}

\begin{proof}
    Set $\delta Z_t=Z^1_t-Z^2_t$ and  
    \begin{align*}
        \delta{f}_t=f^1(t,Y^1_t,Z^1_t)-f^2(t,Y^2_t,Z^2_t).
    \end{align*}
    For any $\alpha\in(1,2]$ and any $\varepsilon>0$, set $\varepsilon_\alpha=\varepsilon(1-\frac{\alpha}{2})$ and $\tilde{Y}_t=|\delta Y_t|^2+\varepsilon_\alpha$. Applying It\^{o}'s formula to $\tilde{Y}_t^{\alpha/2}e^{\beta t}$, we obtain that 
\begin{align*}
    &\tilde{Y}_t^{\frac{\alpha}{2}}e^{\beta t}+\beta\int_t^\tau e^{\beta s}\tilde{Y}^{\frac{\alpha}{2}}_s ds+\frac{\alpha}{2}\int_t^\tau e^{\beta s}\tilde{Y}^{\frac{\alpha}{2}-1}_s|\delta Z_s|^2ds\\
    =&(|\delta \xi|^2+\varepsilon_\alpha)^{\frac{\alpha}{2}}e^{\beta \tau}+\alpha(1-\frac{\alpha}{2})\int_t^\tau e^{\beta s}\tilde{Y}^{\frac{\alpha}{2}-2}|\delta Y_s|^2|\delta Z_s|^2ds\\
    &+\alpha\int_t^\tau e^{\beta s}\tilde{Y}^{\frac{\alpha}{2}-1}_s\delta Y_s\delta f_sds-\alpha\int_t^\tau e^{\beta s}\tilde{Y}^{\frac{\alpha}{2}-1}_s\delta Y_s\delta Z_sdB_s.
    \end{align*}
    Since $f^1$ is Lipschitz continuous, by H\"{o}lder's inequality and Young's inequality, we have 
    \begin{align*}
        &\alpha\int_t^\tau e^{\beta s}\tilde{Y}^{\frac{\alpha}{2}-1}_s\delta Y_s\delta f_sds\leq \alpha\int_t^\tau e^{\beta s}\tilde{Y}^{\frac{\alpha-1}{2}}_s|\delta f_s|ds\\
        \leq &\alpha\int_t^\tau e^{\beta s}\tilde{Y}^{\frac{\alpha-1}{2}}_s(|\delta_2 f_s|+C|\delta Y_s|+C|\delta Z_s|)ds\\
        \leq &(\alpha C+\frac{\alpha C^2}{\alpha-1})\int_t^\tau e^{\beta s}\tilde{Y}^{\frac{\alpha}{2}}_sds+\frac{\alpha(\alpha-1)}{4}\int_t^\tau e^{\beta s}\tilde{Y}^{\frac{\alpha}{2}-1}_s|\delta Z_s|^2ds\\
        &+(\alpha-1)\int_t^\tau e^{\beta s}\tilde{Y}^{\frac{\alpha}{2}}_sds+\int_t^\tau e^{\beta s}|\delta_2 f_s|^\alpha ds.
    \end{align*}
    Therefore, as long as $\beta>\alpha C+\alpha-1+\frac{\alpha C^2}{\alpha-1}$, all the above analysis implies that 
    \begin{align*}
        \tilde{Y}_t^{\frac{\alpha}{2}}e^{\beta t}\leq (|\delta \xi|^2+\varepsilon_\alpha)^{\frac{\alpha}{2}}e^{\beta \tau}+\int_t^\tau e^{\beta s}|\delta_2 f_s|^\alpha ds-\alpha\int_t^\tau e^{\beta s}\tilde{Y}^{\frac{\alpha}{2}-1}_s\delta Y_s\delta Z_sdB_s.
    \end{align*}
    Taking conditional expectations on both sides and letting $\varepsilon$ tend to $0$, we obtain the desired result.
\end{proof}

Now, we introduce the conditional expectation BSDEs. Given a subfiltration $\mathbb{G}=\{\mathcal{G}_t\}_{t\in[0,T]}$, that is, $\mathcal{G}_t\subset \mathcal{F}_t$ for any $t\in[0,T]$ and $\mathcal{G}_s\subset\mathcal{G}_t$ for any $0\leq s\leq t\leq T$, the conditional expectation BSDE with terminal value $\xi$ and generator $f$ takes the following form

\begin{equation}\label{CEBSDE}
    Y_t=\xi+\int_t^T f(s,Y_s,Z_s,\E[Y_s|\mathcal{G}_s],\E[Z_s|\mathcal{G}_s])ds-\int_t^T Z_s dB_s,
\end{equation}
where $\E[\cdot|\mathcal{G}_s]$ represents the conditional expectation w.r.t. $\mathcal{G}_s$ under $\P$.

\begin{assumption}\label{assf}
    The generator $f$ is a map from $\Omega\times[0,T]\times \mathbb{R}\times\mathbb{R}\times\mathbb{R}\times\mathbb{R}$ to $\mathbb{R}$ such that 
    \begin{itemize}
        \item for 
each fixed $(y,z,y',z')$, $f(\cdot,\cdot,y,z,y',z')$ is progressively measurable and $\{f(t,0,0,0,0)\}_{t\in[0,T]}\in\mathcal{H}^2$,
\item there exists a constant $\lambda>0$ such that for any $t\in[0,T]$ and any $y_i,y'_i,z_i,z'_i\in\mathbb{R}$, $i=1,2$
\begin{align*}
|f(t,y_1,z_1,y'_1,z'_1)-f(t,y_2,z_2,y'_2,z'_2)|\leq \lambda(|y_1-y_2|+|z_1-z_2|+|y'_1-y'_2|+|z'_1-z'_2|).
\end{align*}
    \end{itemize} 
\end{assumption}

\begin{theorem}(\cite{Li25})\label{wellposedness for conditional BSDE}
    Given $\xi\in L^2(\mathcal{F}_T)$, suppose that the generator $f$ satisfies Assumption \ref{assf}. Then, the conditional expectation BSDE \eqref{CEBSDE} admits a unique solution $(Y,Z)\in \mathcal{S}^2_c\times\mathcal{H}^2$. 
\end{theorem}

The comparison theorem for conditional expectation BSDEs also holds.

\begin{theorem}(\cite{Li25})\label{comparison for conditional BSDE}
    Let $f^i$, $i=1,2$, be two generators satisfying Assumption \ref{assf}. Moreover, we assume that 
    \begin{itemize}
        \item[(i)] One of the generators is independent of $z'$;
        \item[(ii)] One of the generators is nondecreasing in $y'$.
    \end{itemize}
    Given $\xi^i\in L^2(\mathcal{F}_T)$, $i=1,2$, let $(Y^i,Z^i)$ be the solution to the following conditional expectation BSDE 
    \begin{equation*}\begin{split}\label{CEBSDEK}
        Y^i_t=\xi^i+\int_t^T f^i(s,Y^i_s,Z^i_s,\E[Y^i_s|\mathcal{G}_s],\E[Z^i_s|\mathcal{G}_s])ds
        -\int_t^T Z^i_s dB_s+K^i_T-K^i_t,
    \end{split}\end{equation*}
    where $K^i\in \mathcal{S}^2$, $i=1,2$. Suppose that $\xi^1\leq \xi^2$, $\{K^2_t-K^1_t\}_{t\in[0,T]}$ is nondecreasing and for any $(\omega,t,y,z,y',z')\in\Omega\times[0,T]\times\mathbb{R}\times\mathbb{R}\times\mathbb{R}\times\mathbb{R}$, we have
   \begin{align*}
       f^1(\omega,t,y,z,y',z')\leq f^2(\omega,t,y,z,y',z').
   \end{align*} Then, for all $t\in[0,T]$, we have $Y^1_t\leq Y^2_t$,  $\P$-a.s.
\end{theorem}

\section{Conditional reflected BSDE with Lipschitz generator}

In this section, we consider the conditional reflected BSDE taking the following form:
\begin{equation}\label{conditional reflection}
\begin{cases}
Y_t=\xi+\int_t^T f(s,Y_s,\E[Y_s|\mathcal{G}_s],Z_s)ds-\int_t^T Z_s dB_s+K_T-K_t, \\
Y_t\geq h(t,Y_t,\E[Y_t|\mathcal{G}_t]), \ t\in[0,T],\\
\int_0^T \left(Y_t- h(t,Y_t,\E[Y_t|\mathcal{G}_t])\right)dK_t=0.
\end{cases}
\end{equation}
We propose the following conditions for the sub-filtration $\mathbb{G}=\{\mathcal{G}_t\}_{t\in[0,T]}$, the generator $f$,  the loss function $h$ and the terminal value $\xi$.

\begin{assumption}\label{assG}
    $\mathbb{G}$ is a sub-filtration of $\mathbb{F}$ such that 
    \begin{itemize}
        \item[(a)] $\mathbb{G}$ is nondecreasing and right-continuous,
        \item [(b)] $\mathbb{G}$ is left-quasi-continuous, i.e., $\mathbb{G}$ is left-continuous along stopping times.
    \end{itemize}
\end{assumption}

\begin{assumption}\label{assGfh}
\begin{itemize}
    \item[(i)]The generator $f$ is a map from $\Omega\times[0,T]\times \mathbb{R}\times\mathbb{R}\times\mathbb{R}$ to $\mathbb{R}$ such that 
\begin{itemize}
    \item [(a)] for any $y,y',z\in\mathbb{R}$, $\{f(t,y,y',z)\}_{t\in[0,T]}$ is progressively measurable and $\{f(t,0,0,0)\}_{t\in[0,T]}\in\mathcal{H}^2$,
    \item [(b)] there exists a constant $L>0$ such that $\P$-a.s. for any $t\in[0,T]$ and any $y_i, y'_i,z_i\in\mathbb{R}$, $i=1,2$
\begin{align*}
|f(t,y_1,z_1)-f(t,y_2,z_2)|\leq L(|y_1-y_2|+|y'_1-y'_2|+|z_1-z_2|).
\end{align*}
\end{itemize}
\item[(ii)] $h$ is a mapping from $[0,T]\times\Omega\times\mathbb{R}\times \mathbb{R}$ into $\mathbb{R}$ such that 
    \begin{itemize}
        \item[(a)] for any $y,y'\in\mathbb{R}$, the process $\{h(t,y,y')\}_{t\in[0,T]}$ belongs to $\mathcal{S}^2_c$,
        \item[(b)] there exists two positive constants $\gamma_1,\gamma_2$ such that $\P$-a.s. for all $t\in[0,T]$ and $y_1,y_2,y'_1,y'_2\in\mathbb{R}$, 
        \begin{align*}
            |h(t,y_1,y'_1)-h(t,y_2,y'_2)|\leq \gamma_1|y_1-y_2|+\gamma_2|y'_1-y'_2|.
        \end{align*}
    \end{itemize}
    \item[(iii)] $\xi\in L^2(\mathcal{F}_T)$ and it satisfies $\xi\geq h(T,\xi,\E[\xi|\mathcal{G}_T])$.
\end{itemize}
\end{assumption}

\begin{remark}\label{rmk1}
(i) In the rest of this paper, unless stated, we always suppose that Assumption \ref{assG} holds.

\noindent (ii) Given $S\in\mathcal{S}^2$, the optional projection of $S$ with respect to $\mathbb{G}$ (i.e., the right-continuous modification of $t\rightarrow\E[S_t|\mathcal{G}_t]$) exists (see Chaperter VI in \cite{DM}). In the rest of this paper,  $\{\E[S_t|\mathcal{G}_t]\}_{t\in[0,T]}$ always represents the optional projection. Applying Doob's inequality, we have 
\begin{align*}
    \E\left[\sup_{t\in[0,T]}|\E[S_t|\mathcal{G}_t]|^2\right]\leq\E\left[\sup_{t\in[0,T]}\left|\E\left[\sup_{s\in[0,T]}|S_s|\Big|\mathcal{G}_t\right]\right|^2\right] \leq 4\E\left[\sup_{t\in[0,T]}|S_t|^2\right]<\infty.
\end{align*}
As claimed in \cite{HHL}, when $S\in\mathcal{S}^2_c$  and $\mathbb{G}$ satisfy Assumption \ref{assG},  $\{\E[S_t|\mathcal{G}_t]\}_{t\in[0,T]}$ is also continuous. 
Therefore, if $S\in\mathcal{S}^2_c$, we have $\{\E[S_t|\mathcal{G}_t]\}_{t\in[0,T]}\in\mathcal{S}^2_c$.
\end{remark}

\begin{remark}\label{example}
  Suppose that $\mathbb{G}$ degenerates into the deterministic scenario, i.e., for any $t\in[0,T]$, $\mathcal{G}_t=\mathcal{F}_0$. Then, our conditional reflected BSDE degenerates into the mean-field reflected BSDE in \cite{DEH}, which is closely related with the one in \cite{DDZ,HMW}. 
  
  As shown in \cite{DD,DEH}, the mean-field reflected BSDE is motivated by the applications in pricing life insurance contracts with surrender options. More precisely, consider a portfolio of a large number $N$ of homogeneous life insurance policies $\ell$. Let $(Y^{\ell,N},Z^{\ell,N})$ be the characteristics of the prospective reserve of each policy $\ell=1,2,\cdots,N$. Life insurance constitutes a financial mechanism that embodies the cooperative nature inherent in an insurance contract pool. In this context, prospective reserves are established and priced according to the averaging principle, under which an individual reserve $Y^{\ell,N}$ is benchmarked against the average reserve $\frac{1}{N}\sum_{j=1}^N Y^{j,N}$, commonly referred to as a model point within actuarial practice. Consequently, this study adopts a nonlinear reserve framework, in which the generator $f$ (defined as the reward per unit time) and the solvency level $h$ are jointly determined by the reserve of the specific contract and the average reserve characteristics across $N$ contracts. Given that $N$ is sufficiently large, averaging over the remaining $N-1$ policies yields an effect nearly identical to averaging over all $N$ policies. 
  More specifically, for any $\ell=1,2,\cdots,N$, we have
  \begin{equation}\label{fh}\begin{split}
      f(t,Y^{\ell,N}_t,(Y^{m,N}_t)_{m\neq \ell})&:=\alpha_t-\delta_t Y_t+\beta_t \max \left(\theta_t,Y^{\ell,N}_t-\frac{1}{N}\sum_{j=1}^N Y^{j,N}_t\right),\\
      h(Y^{\ell,N}_t,(Y^{m,N}_t)_{m\neq \ell})&:=u-c(Y^{\ell,N}_t)+\mu\left(\frac{1}{N}\sum_{j=1}^N Y^{j,N}_t-u\right)^+,
  \end{split}\end{equation}
  where $\mu\in(0,1)$. Here, $\delta$ is the discount rate. $\alpha,\beta,\theta$ are positive deterministic functions that constitute the elements of the withdrawal option. $u$ is a minimum required for a benchmark return, $c(Y^{\ell,N}_t)$ represents a reserve dependent management fee. $(\frac{1}{N}\sum_{j=1}^N Y^{j,N}_t-u)^+$ is regarded as a bonus option, which is the possible surplus realized by the average of all the contracts involved. 

Letting $N$ go to infinity in \eqref{fh}, we obtain the following forms of the generator and the obstacle as in \cite{DEH}:
\begin{align*}
    f(t,Y_t,\E[Y_t])&:=\alpha_t-\delta_t Y_t+\beta_t \max(\theta_t,Y_t-\E[Y_t]),\\
    h(Y_t,\E[Y_t])&:=u-c(Y_t)+\mu(\E[Y_t]-u)^+.
\end{align*}
Nevertheless, a key implication of the foregoing analysis is that the sequence $\{Y^{\ell,N}_t\}_{\ell=1}^N$ is independent and identically distributed (i.i.d.). 
In typical situations (see \cite{Smyth}),  $\{Y^{\ell,N}_t\}_{\ell=1}^N$ is only a sequence of conditionally independent  random variables with the same conditional distribution (w.r.t. a sub-$\sigma$-field $\mathcal{G}_t$). Then, by the conditional version of Kolmogorov's strong law of large numbers (see Theorem 4.2 in \cite{MNZ}), we have 
\begin{align*}
    \lim_{N\rightarrow\infty} \frac{1}{N}\sum_{j=1}^N Y^{j,N}_t=\E[Y_t|\mathcal{G}_t].
\end{align*}
In this case, the generator and the obstacle are given by
\begin{align*}
    f(t,Y_t,\E[Y_t|\mathcal{G}_t])&:=\alpha_t-\delta_t Y_t+\beta_t \max(\theta_t,Y_t-\E[Y_t|\mathcal{G}_t]),\\
    h(Y_t,\E[Y_t|\mathcal{G}_t])&:=u-c(Y_t)+\mu(\E[Y_t|\mathcal{G}_t]-u)^+.
\end{align*}
\end{remark}

The main objective of this section is to establish the existence and uniqueness result for the solutions to conditional reflected BSDE \eqref{conditional reflection}. For this purpose,  for any $Y\in \mathcal{S}^2_{\beta,c}$ and $t\in[0,T]$, we define 
\begin{align}\label{Phi}
    \Phi(Y)_t:=\esssup_{\tau\in\mathcal{T}_{t,T}}\mathcal{E}^{f\circ Y}_{t,\tau}\left[h(\tau,Y_\tau,\E[Y_t|\mathcal{G}_t]|_{t=\tau})I_{\{\tau<T\}}+\xi I_{\{\tau=T\}}\right],
\end{align}
where $(f\circ Y)(t,y,z)=f(t,y,\E[Y_t|\mathcal{G}_t],z)$. 
We claim that $\Phi(Y)\in \mathcal{S}^2_{\beta,c}$. It follows from Assumption \ref{assGfh} and Remark \ref{rmk1} that $\{h(t,Y_t,\E[Y_t|\mathcal{G}_t])\}_{t\in[0,T]}\in \mathcal{S}^2_c$. Therefore, the reflected BSDE with terminal value $\xi$, generator $f\circ Y$ and lower obstacle $\{h(t,Y_t,\E[Y_t|\mathcal{G}_t])\}_{t\in[0,T]}$ admits a unique solution $(\hat{Y},\hat{Z},\hat{K})\in\mathcal{S}^2_c\times \mathcal{H}^2\times \mathcal{S}^2_{ci}$ (see \cite{KKPPQ}). By classical results on the relation between the first component of the solution to a reflected BSDE and the value function of an optimal stopping problem under nonlinear expectations (see, e.g., Lemma 2.2 in  \cite{HMW}), for any $t\in[0,T]$, we have $\Phi(Y)_t=\hat{Y}_t$, that is, $\Phi(Y)\in\mathcal{S}^2_c=\mathcal{S}^2_{\beta,c}$.

Now, we present the main result in this section.
\begin{theorem}\label{main}
  Suppose that Assumption \ref{assGfh} holds.  If $\gamma_1,\gamma_2$ satisfy that 
  \begin{align}\label{gamma1gamma2'}
      16\gamma_1^2+64\gamma_2^2<1,
  \end{align}
  Then, the conditional reflected BSDE \eqref{conditional reflection} admits a unique solution $(Y,Z,K)\in\mathcal{S}^2_c\times\mathcal{H}^2\times\mathcal{S}^2_{ci}$.
\end{theorem}

\begin{proof}
Under Condition \eqref{gamma1gamma2'}, there exists a constant $\alpha\in(1,2)$, such that 
\begin{align*}
    2^{\frac{2}{\alpha}}\left(\left(\frac{2}{2-\alpha}\right)^{\frac{2}{\alpha}}\gamma_1^2+\left(\frac{2}{2-\alpha}\right)^{\frac{4}{\alpha}}\gamma_2^2\right)<1.
\end{align*}
Let $\beta,\delta>0$ be such that 
\begin{equation}\begin{split}\label{betadelta}
    &\beta>\alpha L+\alpha-1+\frac{\alpha L^2}{\alpha-1},\\
    &2^{\frac{2}{\alpha}}\gamma_1^2\left(\frac{2}{2-\alpha}\right)^{\frac{2}{\alpha}}+\left(2^{\frac{2}{\alpha}}\gamma_2^2+2^{\frac{2}{\alpha}-1}L^2\delta^{\frac{2}{\alpha}}\right)\left(\frac{2}{2-\alpha}\right)^{\frac{4}{\alpha}}<1.
\end{split}\end{equation}

\textbf{Step 1.} We show that $\Phi:\mathcal{S}^2_{{\frac{2\beta}{\alpha}},c}(T-\delta,T)\rightarrow\mathcal{S}^2_{{\frac{2\beta}{\alpha}},c}(T-\delta,T)$ is a contraction on the time interval $[T-\delta,T]$.  Given $Y,\bar{Y}\in\mathcal{S}^2_{{\frac{2\beta}{\alpha}},c}(T-\delta,T)$, for any $t\in[T-\delta,T]$, by Proposition \ref{a priori estimates}, we have
\begin{equation*}
    \begin{split}
        &e^{\beta t}|\Phi(Y)_t-\Phi(\bar{Y})_t|^\alpha\\
        \leq &e^{\beta t}\esssup_{\tau\in\mathcal{T}_{t,T}}\Big|\mathcal{E}^{f\circ Y}_{t,\tau}\left[h(\tau,Y_\tau,\E[Y_t|\mathcal{G}_t]|_{t=\tau})I_{\{\tau<T\}}+\xi I_{\{\tau=T\}}\right]\\
         &\ \ \ \ \ \ \ \ \ \ \ \ \ \ \ \  -\mathcal{E}^{f\circ \bar{Y}}_{t,\tau}\left[h(\tau,\bar{Y}_\tau,\E[\bar{Y}_t|\mathcal{G}_t]|_{t=\tau})I_{\{\tau<T\}}+\xi I_{\{\tau=T\}}\right]\Big|^\alpha\\
        \leq & \esssup_{\tau\in\mathcal{T}_{t,T}}\E\left[e^{\beta\tau}|h(\tau,{Y}_\tau,\E[{Y}_t|\mathcal{G}_t]|_{t=\tau})-h(\tau,\bar{Y}_\tau,\E[\bar{Y}_t|\mathcal{G}_t]|_{t=\tau})|^\alpha\big|\mathcal{F}_t\right]\\
        &+\esssup_{\tau\in\mathcal{T}_{t,T}}\E\left[\int_t^\tau e^{\beta s}|(f\circ Y)(s,Y^\tau_s,Z^\tau_s)-(f\circ\bar{Y})(s,Y^\tau_s,Z^\tau_s)|^\alpha ds\big|\mathcal{F}_t\right],
    \end{split}
\end{equation*}
where $(Y^\tau,Z^\tau)$ is the solution to the BSDE with terminal time $\tau$, terminal value $h(\tau,Y_\tau,\E[Y_t|\mathcal{G}_t]|_{t=\tau})I_{\{\tau<T\}}+\xi I_{\{\tau=T\}}$ and generator $f\circ Y$. 
Simple calculation yields that 
\begin{equation*}
    \begin{split}
        & \esssup_{\tau\in\mathcal{T}_{t,T}}\E\left[e^{\beta\tau}|h(\tau,{Y}_\tau,\E[{Y}_t|\mathcal{G}_t]|_{t=\tau})-h(\tau,\bar{Y}_\tau,\E[\bar{Y}_t|\mathcal{G}_t]|_{t=\tau})|^\alpha\big|\mathcal{F}_t\right]\\
         \leq & 2^{\alpha-1}\esssup_{\tau\in\mathcal{T}_{t,T}}\E\left[e^{\beta\tau}\left(\gamma_1^\alpha|Y_\tau-\bar{Y}_\tau|^\alpha+\gamma_2^\alpha|\E[Y_t|\mathcal{G}_t]-\E[\bar{Y}_t|\mathcal{G}_t]|^\alpha|_{t=\tau}\right)\big|\mathcal{F}_t\right]\\
         \leq & 2^{\alpha-1}\esssup_{\tau\in\mathcal{T}_{t,T}}\E\left[e^{\beta \tau}\gamma_1^\alpha|Y_\tau-\bar{Y}_\tau|^\alpha+\gamma_2^\alpha\E[e^{\beta t}|Y_t-\bar{Y}_t|^\alpha|\mathcal{G}_t]|_{t=\tau}\big|\mathcal{F}_t\right]\\
         \leq & 2^{\alpha-1}\E\left[\sup_{s\in[T-\delta,T]}e^{\beta s}\gamma_1^\alpha|Y_s-\bar{Y}_s|^\alpha+\sup_{s\in[T-\delta,T]}\gamma_2^\alpha\E\left[e^{\beta s}|Y_s-\bar{Y}_s|^\alpha|\mathcal{G}_s\right]\Big|\mathcal{F}_t\right]\\
         \leq & 2^{\alpha-1}\E\left[\sup_{s\in[T-\delta,T]}e^{\beta s}\gamma_1^\alpha|Y_s-\bar{Y}_s|^\alpha+\sup_{s\in[T-\delta,T]}\gamma_2^\alpha\E\left[\sup_{u\in[T-\delta,T]}e^{\beta u}|Y_u-\bar{Y}_u|^\alpha\Big|\mathcal{G}_s\right]\Bigg|\mathcal{F}_t\right]
    \end{split}
\end{equation*}
and 
\begin{equation*}
    \begin{split}
        &\esssup_{\tau\in\mathcal{T}_{t,T}}\E\left[\int_t^\tau e^{\beta s}|(f\circ Y)(s,Y^\tau_s,Z^\tau_s)-(f\circ\bar{Y})(s,Y^\tau_s,Z^\tau_s)|^\alpha ds\big|\mathcal{F}_t\right]\\
        \leq &\esssup_{\tau\in\mathcal{T}_{t,T}}\E\left[\int_t^\tau e^{\beta s}L^\alpha\left(\E\left[|Y_s-\bar{Y}_s|\big|\mathcal{G}_s\right]\right)^\alpha\Big|\mathcal{F}_t\right]\\
        \leq &\delta L^\alpha \E\left[\sup_{s\in[T-\delta,T]}\E\left[e^{\beta s}|Y_s-\bar{Y}_s|^\alpha\big|\mathcal{G}_s\right]\Big|\mathcal{F}_t\right]\\
        \leq &\delta L^\alpha \E\left[\sup_{s\in[T-\delta,T]}\E\left[\sup_{u\in[T-\delta,T]}e^{\beta u}|Y_u-\bar{Y}_u|^\alpha\Big|\mathcal{G}_s\right]\Bigg|\mathcal{F}_t\right].
    \end{split}
\end{equation*}
Set 
\begin{align*}
    \delta_{\alpha,\beta} Y_t= e^{\beta t}|Y_t-\bar{Y}_t|^\alpha.
\end{align*}
It is easy to check that  
\begin{align*}
    &\E\left[\sup_{t\in[T-\delta,T]}e^{\frac{2\beta t}{\alpha}}|\Phi(Y)_t-\Phi(\bar{Y})_t|^2\right]\\
    \leq &2\E\left[\sup_{t\in[T-\delta,T]}\left(\E\left[\gamma_1^\alpha\sup_{s\in[T-\delta,T]}\delta_{\alpha,\beta}Y_s+\gamma^\alpha_2\sup_{s\in[T-\delta,T]}\E\left[\sup_{u\in[T-\delta,T]}\delta_{\alpha,\beta}Y_u\Big|\mathcal{G}_s\right]\right]\Bigg|\mathcal{F}_t\right)^{\frac{2}{\alpha}}\right]\\
    &+ 2^{\frac{2}{\alpha}-1}L^2\delta^{\frac{2}{\alpha}}\E\left[\sup_{t\in[T-\delta,T]}\left(\E\left[\sup_{s\in[T-\delta,T]}\E\left[\sup_{u\in[T-\delta,T]}\delta_{\alpha,\beta}Y_u\big|\mathcal{G}_s\right]\Bigg|\mathcal{F}_t\right]\right)^{\frac{2}{\alpha}}\right]\\
    \leq&2^{\frac{2}{\alpha}}\gamma_1^2\E\left[\sup_{t\in[T-\delta,T]}\left(\E\left[\sup_{s\in[T-\delta,T]}\delta_{\alpha,\beta}Y_s\Big|\mathcal{F}_t\right]\right)^{\frac{2}{\alpha}}\right]\\
    &+\left(2^{\frac{2}{\alpha}}\gamma_2^2+2^{\frac{2}{\alpha}-1}L^2\delta^{\frac{2}{\alpha}}\right)\E\left[\sup_{t\in[T-\delta,T]}\left(\E\left[\sup_{s\in[T-\delta,T]}\E\left[\sup_{u\in[T-\delta,T]}\delta_{\alpha,\beta}Y_u\Big|\mathcal{G}_s\right]\Bigg|\mathcal{F}_t\right]\right)^{\frac{2}{\alpha}}\right].
\end{align*}
Applying Doob's inequality implies that 
\begin{align*}
    \E\left[\sup_{t\in[T-\delta,T]}\left(\E\left[\sup_{s\in[T-\delta,T]}\delta_{\alpha,\beta}Y_s\Bigg|\mathcal{F}_t\right]\right)^{\frac{2}{\alpha}}\right]&\leq \left(\frac{2}{2-\alpha}\right)^{\frac{2}{\alpha}}\E\left[\sup_{t\in[T-\delta,T]}\delta_{\frac{2\beta}{\alpha},2}Y_t\right]\\
    &=\left(\frac{2}{2-\alpha}\right)^{\frac{2}{\alpha}}\E\left[\sup_{t\in[T-\delta,T]}e^{\frac{2\beta t}{\alpha}}|Y_t-\bar{Y}_t|^2\right] 
\end{align*}
and 
\begin{align*}
    &\E\left[\sup_{t\in[T-\delta,T]}\left(\E\left[\sup_{s\in[T-\delta,T]}\E\left[\sup_{u\in[T-\delta,T]}\delta_{\alpha,\beta}Y_u\Bigg|\mathcal{G}_s\right]\Bigg|\mathcal{F}_t\right]\right)^{\frac{2}{\alpha}}\right]\\
    \leq& \left(\frac{2}{2-\alpha}\right)^{\frac{2}{\alpha}}\E\left[\sup_{s\in[T-\delta,T]}\left(\E\left[\sup_{u\in[T-\delta,T]}\delta_{\alpha,\beta}Y_u\Bigg|\mathcal{G}_s\right]\right)^{\frac{2}{\alpha}}\right]\\
    \leq &\left(\frac{2}{2-\alpha}\right)^{\frac{4}{\alpha}}\E\left[\sup_{t\in[T-\delta,T]}\delta_{\frac{2\beta}{\alpha},2}Y_t\right]
    =\left(\frac{2}{2-\alpha}\right)^{\frac{4}{\alpha}}\E\left[\sup_{t\in[T-\delta,T]}e^{\frac{2\beta t}{\alpha}}|Y_t-\bar{Y}_t|^2\right]. 
\end{align*}
All the above analysis indicates that 
\begin{align*}
    \E\left[\sup_{t\in[T-\delta,T]}e^{\frac{2\beta t}{\alpha}}|\Phi(Y)_t-\Phi(\bar{Y})_t|^2\right]\leq C(\gamma_1,\gamma_2,L,\delta)\E\left[\sup_{t\in[T-\delta,T]}e^{\frac{2\beta t}{\alpha}}|Y_t-\bar{Y}_t|^2\right],
\end{align*}
where 
\begin{align*}
    C(\gamma_1,\gamma_2,L,\delta)=2^{\frac{2}{\alpha}}\gamma_1^2\left(\frac{2}{2-\alpha}\right)^{\frac{2}{\alpha}}+\left(2^{\frac{2}{\alpha}}\gamma_2^2+2^{\frac{2}{\alpha}-1}L^2\delta^{\frac{2}{\alpha}}\right)\left(\frac{2}{2-\alpha}\right)^{\frac{4}{\alpha}}.
\end{align*}
Hence, for $\delta$ satisfying \eqref{betadelta}, $\Phi$ is a contraction. 

\textbf{Step 2.} We show that the conditional reflected BSDE \eqref{conditional reflection} has a unique solution on $[T-\delta,T]$. By Step 1, there exists a unique $Y^0\in \mathcal{S}^2_{{\frac{2\beta}{\alpha}},c}(T-\delta,T)=\mathcal{S}^2_c(T-\delta,T)$ such that, for any $t\in[T-\delta,T]$
\begin{align*}
    Y^0_t=\esssup_{\tau\in\mathcal{T}_{t,T}}\mathcal{E}^{f\circ Y^0}_{t,\tau}\left[h(\tau,Y^0_\tau,\E[Y^0_t|\mathcal{G}_t]|_{t=\tau})I_{\{\tau<T\}}+\xi I_{\{\tau=T\}}\right].
\end{align*}
Now, let $(\hat{Y},\hat{Z},\hat{K})\in\mathcal{S}^2_c\times\mathcal{H}^2\times\mathcal{S}^2_{ci}$ be the unique solution to the reflected BSDE with terminal value $\xi$, generator $g$ and obstacle $\{h(t,Y^0_t,\E[Y^0_t|\mathcal{G}_t])\}_{t\in[0,T]}$ with $g(t,y,z)=f(t,y,\E[Y^0_t|\mathcal{G}_t],z)$ on $[T-\delta,T]$. Using the link between reflected BSDEs and optimal stopping problem under nonlinear expectations, we have 
\begin{align*}
    \hat{Y}_t=\esssup_{\tau\in\mathcal{T}_{t,T}}\mathcal{E}^{g}_{t,\tau}\left[h(\tau,Y^0_\tau,\E[Y^0_t|\mathcal{G}_t]|_{t=\tau})I_{\{\tau<T\}}+\xi I_{\{\tau=T\}}\right].
\end{align*}
It follows that $Y^0\equiv \hat{Y}$. Set $Z^0=\hat{Z}$ and $K^0=\hat{K}$. Clearly, $(Y^0,Z^0,K^0)$ is a solution to conditional reflected BSDE \eqref{conditional reflection} on $[T-\delta,T]$. Since $Y^0$ is a nonlinear $g$-supermartingale, the uniqueness of the nonlinear Doob-Meyer decomposition (see Theorem 3.3 in \cite{P99}) yields the uniqueness of $(Z^0,K^0)$. 

\textbf{Step 3.} We show the existence and uniqueness on the whole time interval $[0,T]$. Set $N=[\frac{T}{\delta}]$. Using the same reasoning on each time interval $[(T-(i+1)\delta)\vee 0,T-i\delta]$ with a similar dynamics but terminal condition $Y^{i-1}_{T-i\delta}$ at time $T-i\delta$, we build recursively for $i=1$ to $N$ a solution $(Y^i,Z^i,K^i)$ on each time interval $[(T-(i+1)\delta)\vee 0,T-i\delta]$. Pasting properly these processes, we naturally derive a solution $(Y,Z,K)$ satisfying \eqref{conditional reflection} on
the whole time interval $[0,T ]$. The uniqueness of the solution on the whole interval is a direct consequence of the one on each small interval. The proof is complete. 
\end{proof}

One of the important properties of the classical reflected BSDEs and BSDEs with conditional reflection is the comparison theorem (see, e.g., Theorem 4.1 in \cite{KKPPQ} and Corollary 3.3 in \cite{HHL}). The comparison theorem for conditional reflected BSDEs still holds. 
\begin{theorem}\label{comparison}
    Suppose that $\xi^i,f^i,h^i$ satisfy Assumption \ref{assGfh}, $i=1,2$.  Let $(Y^i,Z^i,K^i)$ be the solution to the conditional reflected BSDE with parameters $\xi^i,f^i,h^i$. Moreover, suppose that 
    \begin{itemize}
        \item[(1)] one of the generators is non-decreasing in $y'$,
        \item[(2)] one of the loss functions satisfies the following condition
        \begin{itemize}
            \item for any fixed $(t,y)$,  the mapping $y'\rightarrow h^i(t,y,y')$ is non-increasing,
            \item for any fixed $(t,y')$, the mapping $y\rightarrow h^i(t,y,y')$ is non-increasing.
        \end{itemize}
    \end{itemize}
     If $\xi^1\leq \xi^2$, a.s., for any $(y,y',z)\in\mathbb{R}^3$, $f^1(t,y,y',z)\leq f^2(t,y,y',z)$, $d\P\times dt$-a.e. and $h^1(t,y,y')\leq h^2(t,y,y')$, $d\P\times dt$-a.e., we have $Y^1_t\le Y^2_t$, $t\in[0,T]$. 
\end{theorem}

\begin{proof}
   Without loss of generality, suppose that $f^1$ is non-decreasing in $y'$ and $h^1$ satisfy condition (2). Applying It\^{o}'s formula to $|(Y^1_t-Y^2_t)^+|^2$ and taking expectation, we have 
    \begin{align*}
        &\E\left[|(Y^1_t-Y^2_t)^+|^2\right]+\E\left[\int_t^T I_{A_s} |Z^1_s-Z^2_s|^2ds\right]\\
        \leq &2\E\left[\int_t^T(Y^1_s-Y^2_s)^+\delta f_sds\right] +2\E\left[\int_t^T(Y^1_s-Y^2_s)^+(dK^1_s-dK^2_s)\right],
    \end{align*}
    where $A_s=\{Y^1_s>Y^2_s\}$ and
    \begin{align*}
        \delta f_s=f^1(s,Y^1_s,\E[Y^1_s|\mathcal{G}_s],Z^1_s)-f^2(s,Y^2_s,\E[Y^2_s|\mathcal{G}_s],Z^2_s).
    \end{align*}
     Since on $A_t$, we have 
     \begin{align*}
         Y^1_t>Y^2_t\geq h^2(t,Y^2_t,\E[Y^2_t|\mathcal{G}_t])\geq h^1(t,Y^2_t,\E[Y^2_t|\mathcal{G}_t])\geq h^1(t,Y^1_t,\E[Y^1_t|\mathcal{G}_t]).
     \end{align*}
    It follows that 
    \begin{align*}
        &\int_t^T(Y^1_s-Y^2_s)^+(dK^1_s-dK^2_s)\\
        \leq &\int_t^T I_{A_s}(Y^1_s-h^1(s,Y^1_s,\E[Y^1_s|\mathcal{G}_s]))dK^1_s-\int_t^T(Y^1_s-Y^2_s)^+dK^2_s\leq 0.
    \end{align*}
    All the above analysis indicates that 
    \begin{align*}
        &\E\left[|(Y^1_t-Y^2_t)^+|^2\right]+\E\left[\int_t^T I_{A_s} |Z^1_s-Z^2_s|^2ds\right]\\
        \leq &2\E\left[\int_t^T(Y^1_s-Y^2_s)^+(f^1(s,Y^1_s,\E[Y^1_s|\mathcal{G}_s],Z^1_s)-f^1(s,Y^2_s,\E[Y^2_s|\mathcal{G}_s],Z^2_s))ds\right]\\
       \leq &2L\E\left[\int_t^T(Y^1_s-Y^2_s)^+\left(|Y^1_s-Y^2_s|+|Z^1_s-Z^2_s|+\left(\E[Y^1_s-Y^2_s|\mathcal{G}_s]\right)^+\right)ds\right]\\
       \leq &\E\left[\int_t^T I_{A_s} |Z^1_s-Z^2_s|^2ds\right]+(L^2+4L)\E\left[\int_t^T|(Y^1_s-Y^2_s)^+|ds\right],
    \end{align*}
    where we have used the Fubini theorem and the following fact in the last inequality
    \begin{align*}
        \E\left[(Y^1_s-Y^2_s)^+\left(\E[Y^1_s-Y^2_s|\mathcal{G}_s]\right)^+\right]\leq &\E\left[(Y^1_s-Y^2_s)^+\E[(Y^1_s-Y^2_s)^+|\mathcal{G}_s]\right]\\
        \leq& \E\left[\left(\E[(Y^1_s-Y^2_s)^+|\mathcal{G}_s]\right)^2\right]\\
        \leq &\E\left[\left|(Y^1_s-Y^2_s)^+\right|^2\right].
    \end{align*}
    Therefore, we have
    \begin{align*}
        \E\left[|(Y^1_t-Y^2_t)^+|^2\right]\leq (L^2+4L)\E\left[\int_t^T|(Y^1_s-Y^2_s)^+|ds\right].
    \end{align*}
    Applying Gronwall's inequality, we obtain $(Y^1_t-Y^2_t)^+=0$, $t\in[0,T]$. The proof is complete. 
\end{proof}

\begin{remark}
    (i) As far as we know, the comparison theorem for mean-field reflected BSDE with constraint $Y_t\geq h(t,Y_t,\E[Y_t])$ (or $Y_t\geq h(t,Y_t,\P_{Y_t})$) has not been studied yet. Theorem \ref{comparison} is the first attempt to establish the comparison property.

 \noindent   (ii) The comparison theorem for mean-field BSDE needs the assumption that the generator is non-decreasing for $y'$ (see Theorem 3.2 and Example 3.2 in \cite{BLP}). Therefore, condition (1) in Theorem \ref{comparison} is necessary. 

    \noindent (iii) Suppose that $h$ does not depend on $y'$ and  may be increasing in $y$. Note that the comparison theorem does not need condition \eqref{gamma1gamma2'}, the function $y-h(t,y)$ may be decreasing in $y$. In this case, in order to fulfill the constraint $Y_t-h(t,Y_t)\geq 0$, the compensating term $K$ may be decreasing, which is unreasonable. 

     \noindent (iv) For the mean reflected problem, the constraint is given by $\E[l(t,Y_t)]\geq 0$, where  $l$ is a loss function, which is strictly increasing in its last component. If we choose our loss function $h(t,y,y')=y-l(t,y')$, our constraint can be written as 
   \begin{align*}
        l(t,\E[Y_t|\mathcal{G}_t])\geq 0.
    \end{align*}
    Therefore, it is reasonable that the loss function $h$ is non-increasing in $y'$. 
\end{remark}

\section{Conditional reflected BSDE with quadratic generator}

In this section, we consider the conditional reflected BSDE \eqref{conditional reflection} when the generator has quadratic growth in its $z$-component. For the case of bounded terminal value and loss function, we will apply a linearization technique, a fixed-point argument and the BMO martingale theory to obtain the well-posedness result. For the case of unbounded terminal value and loss function, we only show the uniqueness of the solution by the $\theta$-method. 

\subsection{Bounded terminal value and loss function: the existence and uniqueness}

In this subsection, we consider the conditional reflected BSDE \eqref{conditional reflection} with parameters $(\xi,f,h)$, which satisfy the following condition.
\begin{assumption}\label{ass quadratic}
\begin{itemize}
    \item[(i)]The generator $f$ is a map from $\Omega\times[0,T]\times \mathbb{R}\times\mathbb{R}\times\mathbb{R}$ to $\mathbb{R}$ such that 
\begin{itemize}
    \item [(a)] for any $y,y',z\in\mathbb{R}$, $\{f(t,y,y',z)\}_{t\in[0,T]}$ is progressively measurable and there exist positive constants $\alpha,\beta,\gamma$ such that for any $t\in[0,T]$, $y,y',z\in\mathbb{R}$,
    \begin{align*}
        |f(t,y,y',z)|\leq \alpha+\beta(|y|+|y'|)+\frac{\gamma}{2}|z|^2,
    \end{align*}
    \item [(b)] there exists a constant $L>0$ such that $\P$-a.s. for any $t\in[0,T]$ and any $y_i, y'_i,z_i\in\mathbb{R}$, $i=1,2$,
\begin{align*}
|f(t,y_1,y'_1,z_1)-f(t,y_2,y'_2,z_2)|\leq L(|y_1-y_2|+|y'_1-y'_2|)+L(1+|z_1|+|z_2|)|z_1-z_2|.
\end{align*}
\end{itemize}
\item[(ii)] $h$ is a mapping from $[0,T]\times\Omega\times\mathbb{R}\times \mathbb{R}$ into $\mathbb{R}$ such that 
    \begin{itemize}
        \item[(a)] for any $y,y'\in\mathbb{R}$, the process $\{h(t,y,y')\}_{t\in[0,T]}\in\mathcal{S}^\infty_c$ and $h$ is uniformly bounded w.r.t. $(t,\omega,y,y')$,
        \item[(b)] there exists two positive constants $\gamma_1,\gamma_2$ such that $\P$-a.s. for all $t\in[0,T]$ and $y_1,y_2,y'_1,y'_2\in\mathbb{R}$, 
        \begin{align*}
            |h(t,y_1,y'_1)-h(t,y_2,y'_2)|\leq \gamma_1|y_1-y_2|+\gamma_2|y'_1-y'_2|.
        \end{align*}
    \end{itemize}
    \item[(iii)] $\xi\in L^\infty(\mathcal{F}_T)$ and it satisfies $\xi\geq h(T,\xi,\E[\xi|\mathcal{G}_T])$.
\end{itemize}
\end{assumption}

Now, we are ready to state the main result of this subsection.

\begin{theorem}\label{thm quadratic}
    Suppose that Assumption \ref{ass quadratic} holds. If $\gamma_1,\gamma_2$ satisfy 
    \begin{align}\label{gamma1gamma2''}
        \gamma_1+\gamma_2<1,
    \end{align}
    then, the quadratic conditional reflected BSDE \eqref{conditional reflection}  admits a unique solution $(Y,Z,K)\in\mathcal{S}^\infty_c\times BMO\times \mathcal{S}_{ci}$.
\end{theorem}

We first provide some ideas before proving the main result. Given $U\in\mathcal{S}^\infty_c$, for any $(t,z)$, set $(f\circ U)(t,z)=f(t,U_t,\E[U_t|\mathcal{G}_t],z)$ and $(h\circ U)_t=h(t,U_t,\E[U_t|\mathcal{G}_t])$. Consider the quadratic reflected BSDE with terminal value $\xi$, generator $f\circ U$ and obstacle $h\circ U$. By a similar analysis as the proof of Lemma 3.3 in \cite{HMW}, it has a unique solution $(Y,Z,K)\in \mathcal{S}^\infty_c\times BMO\times \mathcal{S}_{ci}$. By Lemma 2.2 in \cite{HMW}  or Proposition 3.1 in \cite{KLQT}, for any $t\in[0,T]$, we have  
\begin{align*}
    Y_t=\esssup_{\tau\in\mathcal{T}_{t,T}}y^{\tau,U}_t(=:\Phi(U)_t),
\end{align*}
where $y^{\tau,U}\in\mathcal{S}^\infty_c$ is the first component of the solution to the quadratic BSDE with terminal value $(h\circ U)_\tau I_{\{\tau<T\}}+\xi I_{\{\tau=T\}}$ and generator $f\circ U$. Now, we have constructed a map $\Phi:\mathcal{S}^\infty_c\rightarrow\mathcal{S}^\infty_c$. Its fixed point is exactly the first component of the solution to the quadratic conditional reflected BSDE \eqref{conditional reflection}.
 
\begin{proof}[Proof of Theorem \ref{thm quadratic}]
    Given $Y,\bar{Y}\in \mathcal{S}^\infty_c$, for any $t\in[0,T]$ and $\tau\in\mathcal{T}_{t,T}$, let $(y^\tau,z^\tau)$ (resp., $(\bar{y}^\tau,\bar{z}^\tau)$) be the solution to the BSDE on the time horizon $[0,\tau]$ with terminal value $\xi I_{\{\tau=T\}}+h(\tau,Y_\tau,\E[Y_t|\mathcal{G}_t]|_{t=\tau})$ (resp. $\xi I_{\{\tau=T\}}+h(\tau,\bar{Y}_\tau,\E[\bar{Y}_t|\mathcal{G}_t]|_{t=\tau})$) and generator $f\circ Y$ (resp. $f\circ {\bar{Y}}$), where $
    (f\circ Y)(t,z)=f(t,Y_t,\E[Y_t|\mathcal{G}_t],z)$ (resp. $(f\circ {\bar{Y}})(t,z)=f(t,\bar{Y}_t,\E[\bar{Y}_t|\mathcal{G}_t],z)$). Similar with \eqref{Phi}, for any $t\in[0,T]$, we define
    \begin{equation*}
        \Phi(Y)_t:=\esssup_{\tau\in\mathcal{T}_{t,T}}y^\tau_t.
    \end{equation*}
    We first prove that $\Phi:\mathcal{S}^\infty_c\rightarrow \mathcal{S}^\infty_c$ is a contraction map on $[T-\delta,T]$ when $\delta$ satisfies 
    \begin{align*}
        \gamma_1+\gamma_2+2L\delta<1.
    \end{align*}

    For any $t\in[0,T]$, set 
    \begin{align*}
        \beta_t=\frac{(f\circ{Y})(t,z^\tau_t)-(f\circ Y)(t,\bar{z}^\tau_t)}{|z^\tau_t-\bar{z}^\tau_t|^2}(z^\tau_t-\bar{z}^\tau_t)I_{\{z^\tau_t\neq \bar{z}^\tau_t\}}.
    \end{align*}
    It is easy to check that $\beta\in BMO$ and $(y^\tau-\bar{y}^\tau,z^\tau-\bar{z}^\tau)$ solves the following BSDE:
    \begin{align*}
        y^\tau_t-\bar{y}^\tau_t=&h(\tau,Y_\tau,\E[Y_t|\mathcal{G}_t])|_{t=\tau}-h(\tau,\bar{Y}_\tau,\E[\bar{Y}_t|\mathcal{G}_t])|_{t=\tau}-\int_t^T I_{[0,\tau]}(s)(z^\tau_s-\bar{z}^\tau_s)dB_s\\
        &+\int_t^T \left(\beta_s(z^\tau_s-\bar{z}^\tau_s)+(f\circ  Y)(s,\bar{z}^\tau_s)-(f\circ{\bar{Y}})(s,\bar{z}^\tau_s)\right)I_{[0,\tau]}(s)ds.
    \end{align*}
    We define $\widetilde{B}_t:=B_t-\int_0^t \beta_s I_{[0,\tau]}(s)ds$ and a probability measure $\widetilde{\P}$ with density given by 
    \begin{align*}
        \frac{d\widetilde{\P}}{d\P}=\mathcal{E}(\beta I_{[0,\tau]}\cdot B)_0^T,
    \end{align*}
    where for any $0\leq s\leq t\leq T$, 
    \begin{align*}
        \mathcal{E}(\beta I_{[0,\tau]}\cdot B)_s^t=\exp\left(\int_s^t \beta_r I_{[0,\tau]}(r)dB_r-\frac{1}{2}\int_s^t \beta_r^2 I_{[0,\tau]}(r)dr\right).
    \end{align*}
    Then, $\widetilde{B}$ is a Brownian motion under $\widetilde{\P}$. Moreover, for any $t\in[T-\delta,T]$ and $\tau\in\mathcal{T}_{t,T}$, we have 
    \begin{align*}
       y^\tau_t-\bar{y}^\tau_t
       =&\E^{\widetilde{\P}}_t\left[h(\tau,Y_\tau,\E[Y_t|\mathcal{G}_t])|_{t=\tau}-h(\tau,\bar{Y}_\tau,\E[\bar{Y}_t|\mathcal{G}_t])|_{t=\tau}\right]\\
       &+\E^{\widetilde{\P}}_t\left[\int_t^T \left((f\circ Y)(s,\bar{z}^\tau_s)-(f\circ{\bar{Y}})(s,\bar{z}^\tau_s)\right)I_{[0,\tau]}(s)ds\right] \\
       \leq &(\gamma_1+L\delta)\|Y-\bar{Y}\|_{\mathcal{S}^\infty(T-\delta,T)}+(\gamma_2+L\delta)\esssup_{t\in[T-\delta,T]}\E[|Y_t-\bar{Y}_t||\mathcal{G}_t]\\
       \leq &(\gamma_1+\gamma_2+2L\delta)\|Y-\bar{Y}\|_{\mathcal{S}^\infty(T-\delta,T)},
    \end{align*}
    where $\E^{\widetilde{\P}}_t[\cdot]$ is the conditional expectation w.r.t. $\mathcal{F}_t$ under $\widetilde{\P}$.    It follows that 
    \begin{align*}
        \|\Phi(Y)-\Phi(\bar{Y})\|_{\mathcal{S}^\infty(T-\delta,T)}\leq (\gamma_1+\gamma_2+2L\delta)\|Y-\bar{Y}\|_{\mathcal{S}^\infty(T-\delta,T)}.
    \end{align*}
    Clearly, $\Phi$ is a contraction map on the time interval $[T-\delta,T]$.

    Now, let $Y\in\mathcal{S}^\infty(T-\delta,T)$ be the fixed point of $\Phi$ on $[T-\delta,T]$, i.e.,
    \begin{align*}
        Y_t=\esssup_{\tau\in\mathcal{T}_{t,T}}\mathcal{E}^{f\circ  Y}_{t,\tau}\left[h(\tau,Y_\tau,\E[Y_t|\mathcal{G}_t]|_{t=\tau})I_{\{\tau<T\}}+\xi I_{\{\tau=T\}}\right].
    \end{align*}
    By Lemma 3.3 in \cite{HMW} (up to a slight modification), the reflected BSDE  with terminal value $\xi$, generator $f\circ Y$ and obstacle $h\circ Y$ admits a unique solution $(\widetilde{Y},Z,K)\in \mathcal{S}^\infty[T-\delta,T]\times  BMO[T-\delta,T]\times\mathcal{A}[T-\delta,T]$. 
    It follows from Lemma 2.2 in \cite{HMW}  or Proposition 3.1 in \cite{KLQT} that $\widetilde{Y}=\Phi(Y)=Y$. Consequently, $(Y,Z,K)$ is a solution to conditional reflected BSDE \eqref{conditional reflection} on the time interval $[T-\delta,T]$. 

    We are now in a position to show the uniqueness. Suppose that $(Y',Z',K')$ is another solution to the conditional reflected BSDE \eqref{conditional reflection} on the time interval $[T-\delta,T]$. Since $Y'$ is the fixed point of $\Phi$, we have $Y'\equiv Y$. Applying It\^{o}'s formula to $|Y_t-Y'_t|^2$, we obtain that $Z\equiv Z'$ and $K\equiv K'$. 

    Finally, a global solution on the entire time interval $[0,T]$ is obtained by stitching the local ones. The global uniqueness can be derived from the local uniqueness on each small time interval. The proof is complete.
\end{proof}

\subsection{Unbounded terminal value and loss function: the uniqueness result}

In this subsection, we consider the conditional reflected BSDE \eqref{conditional reflection} with parameters $(\xi,f,h)$, which satisfy the following condition.
\begin{assumption}\label{ass quadratic'}
\begin{itemize}
    \item[(i)]The generator $f$ is a map from $\Omega\times[0,T]\times \mathbb{R}\times\mathbb{R}\times\mathbb{R}$ to $\mathbb{R}$ such that 
\begin{itemize}
    \item [(a)] for any $y,y',z\in\mathbb{R}$, $\{f(t,y,y',z)\}_{t\in[0,T]}$ is progressively measurable and there exist positive constants $\alpha,\beta,\gamma$ such that for any $t\in[0,T]$, $y,y',z\in\mathbb{R}$,
    \begin{align*}
        |f(t,y,y',z)|\leq \alpha+\beta(|y|+|y'|)+\frac{\gamma}{2}|z|^2,
    \end{align*}
    \item [(b)] there exists a constant $\beta>0$ such that $\P$-a.s. for any $t\in[0,T]$, $y_i, y'_i\in\mathbb{R}$, $i=1,2$ and any $z\in\mathbb{R}$,
\begin{align*}
|f(t,y_1,y'_1,z)-f(t,y_2,y'_2,z)|\leq \beta(|y_1-y_2|+|y'_1-y'_2|),
\end{align*}
\item[(c)] for any $(t,\omega,y,y')\in [0,T]\times\Omega\times\mathbb{R}\times \mathbb{R}$, $f(t,y,y',\cdot)$ is concave or convex.
\end{itemize}
\item[(ii)] $h$ is a mapping from $[0,T]\times\Omega\times\mathbb{R}\times \mathbb{R}$ into $\mathbb{R}$ such that 
    \begin{itemize}
        \item[(a)] for any $y,y'\in\mathbb{R}$, the process $\{h(t,y,y')\}_{t\in[0,T]}\in\mathbb{S}_c$,
        \item[(b)] there exists two positive constants $\gamma_1,\gamma_2$ such that $\P$-a.s. for all $t\in[0,T]$ and $y_1,y_2,y'_1,y'_2\in\mathbb{R}$, 
        \begin{align*}
            |h(t,y_1,y'_1)-h(t,y_2,y'_2)|\leq \gamma_1|y_1-y_2|+\gamma_2|y'_1-y'_2|.
        \end{align*}
    \end{itemize}
    \item[(iii)] $\xi\in \mathbb{L}(\mathcal{F}_T)$ and it satisfies $\xi\geq h(T,\xi,\E[\xi|\mathcal{G}_T])$.
\end{itemize}
\end{assumption}

We first recall the following technical lemma, which plays a key role in the following analysis. 
\begin{lemma}(\cite{FHT,HMW})\label{technical lemma}
   Given $\tau\in\mathcal{T}_{0,T}$, Let $(y^\tau,z^\tau)\in\mathcal{S}^2_c\times\mathcal{H}^2$ be a solution to the following BSDE:
    \begin{align*}
        y^\tau_t=\eta+\int_t^\tau g(s,y^\tau_s,z^\tau_s)ds-\int_t^\tau z^\tau_sdB_s.
    \end{align*}
    Suppose that there is a constant $p\geq 1$ such that 
    \begin{align*}
        \E\left[\exp\left(2p\gamma e^{\beta T}\sup_{t\in[0,T]}|y^\tau_t|+2p\gamma \int_0^T \alpha_t e^{\beta t}dt\right)\right]<\infty,
    \end{align*}
    where $\alpha$ is a progressively measurable process and $\beta,\gamma$ are two nonnegative constants. Then, we have
    \begin{itemize}
        \item[(i)] Suppose that $g$ satisfies the following condition
        \begin{itemize}
            \item[(H1)] for each $(t,\omega,y,z)\in[0,T]\times \Omega\times\mathbb{R}\times\mathbb{R}$, $|g(t,y,z)|\leq \alpha_t+\beta|y|+\frac{\gamma}{2}|z|^2$ with $\alpha\geq 0$ being such that $\int_0^T \alpha_t dt\in \mathbb{L}(\mathcal{F}_T)$ .
        \end{itemize}
        For each $t\in[0,T]$  and $p\geq 1$,
        \begin{align*}
            \exp(p\gamma|y^\tau_t|)\leq \E_t\left[\exp\left(p\gamma e^{\beta( T-t)}|\eta|+p\gamma \int_t^T \alpha_s e^{\beta (s-t)}ds\right)\right].
        \end{align*}
        \item[(ii)] Suppose that $g$ satisfies the following condition
        \begin{itemize}
            \item[(H2)] for each $(t,\omega,y,z)\in[0,T]\times \Omega\times\mathbb{R}\times\mathbb{R}$, $g(t,y,z)\leq \alpha_t+\beta|y|+\frac{\gamma}{2}|z|^2$.
        \end{itemize}
        For each $t\in[0,T]$  and $p\geq 1$,
        \begin{align*}
            \exp(p\gamma(y^\tau_t)^+)\leq \E_t\left[\exp\left(p\gamma e^{\beta( T-t)}\eta^++p\gamma \int_t^T \alpha_s e^{\beta (s-t)}ds\right)\right].
        \end{align*}
    \end{itemize}
\end{lemma}

\begin{lemma}\label{technical lemma2}
    Suppose that $S\in\mathbb{S}_c$. Then, for any $p>1$, we have 
    \begin{align*}
        \E\left[\exp\left(p\sup_{t\in[0,T]}\E\left[|S_t||\mathcal{G}_t\right]\right)\right]\leq \left(\frac{p}{p-1}\right)^{p}\E\left[\exp\left(p\sup_{t\in[0,T]}|S_t|\right)\right].
    \end{align*}
\end{lemma}

\begin{proof}
    Applying Jensen' inequality and Doob's inequality, we obtain that 
\begin{align*}
        \E\left[\exp\left(p\sup_{t\in[0,T]}\E\left[|S_t||\mathcal{G}_t\right]\right)\right]
        \leq &\E\left[\sup_{t\in[0,T]}\left(\exp\left(\E\left[\sup_{s\in[0,T]}|S_s|\Big|\mathcal{G}_t\right]\right)\right)^p\right]\\
        \leq &\E\left[\sup_{t\in[0,T]}\left(\E\left[\exp\left(\sup_{s\in[0,T]}|S_s|\right)\Bigg|\mathcal{G}_t\right]\right)^p\right]\\
        \leq &\left(\frac{p}{p-1}\right)^{p}\E\left[\exp\left(p\sup_{t\in[0,T]}|S_t|\right)\right].
    \end{align*}
\end{proof}

 Given $U\in\mathbb{S}_c$, for any $(t,z)$, set $(f\circ U)(t,z)=f(t,U_t,\E[U_t|\mathcal{G}_t],z)$ and $(h\circ U)_t=h(t,U_t,\E[U_t|\mathcal{G}_t])$. It is easy to check that $h\circ U\in \mathbb{S}_c$ and 
 \begin{itemize}
     \item $f\circ U$ satisfies (H1) in Lemma \ref{technical lemma},
     \item for any $\tau\in \mathcal{T}_{0,T}$, $(h\circ U)_\tau I_{\{\tau<T\}}+\xi I_{\{\tau=T\}}\in \mathbb{L}(\mathcal{F}_T)$.
 \end{itemize}
Therefore, according to Corollary 6 in \cite{BH2}, the quadratic BSDE with terminal value $(h\circ U)_\tau I_{\{\tau<T\}}+\xi I_{\{\tau=T\}}$ and generator $f\circ U$ has a unique solution $(y^{\tau,U},z^{\tau,U})\in\mathbb{S}_c\times\mathcal{H}$. Moreover, as claimed in Remark 4.2 of \cite{HMW}, the comparison theorem for quadratic BSDEs with generator $f\circ U$ also holds. Then, by Lemma 2.2 in \cite{HMW}, for any $t\in[0,T]$, we have  
\begin{align*}
    Y_t=\esssup_{\tau\in\mathcal{T}_{t,T}}y^{\tau,U}_t.
\end{align*}

\begin{proposition}\label{uniqueness}
    Suppose that Assumption \ref{ass quadratic'} holds. If $\gamma_1,\gamma_2$ satisfy
    \begin{equation}\label{eq17'}
        4(\gamma_1+\gamma_2)<1        
    \end{equation}
    the conditional reflected BSDE \eqref{conditional reflection} has at most one solution $(Y,Z,K)\in \mathbb{S}_c\times\mathcal{H}\times\mathcal{S}_{ci}$.
\end{proposition}

\begin{proof}
   Without loss of generality, we assume that $f(t,y,y',\cdot)$ is concave. Let $(Y^i,Z^i,K^i)\in \mathbb{S}_c\times\mathcal{H}\times\mathcal{S}_{ci}$ be the solution to the conditional reflected BSDE \eqref{conditional reflection}, $i=1,2$. Let $(y^{i,\tau},z^{i,\tau})$ be the  solution to BSDE with terminal value $(h\circ Y^i)_\tau I_{\{\tau<T\}}+\xi I_{\{\tau=T\}}$ and generator $f\circ Y^i$. We have 
    \begin{align}\label{eq17}
    Y^i_t=\esssup_{\tau\in\mathcal{T}_{t,T}}y^{i,\tau}_t.
\end{align}

For any $\theta\in(0,1)$, we define
\begin{align*}
    \delta_\theta\ell:=\frac{\theta \ell^1-\ell^2}{1-\theta}, \  \delta_\theta\widetilde{\ell}:=\frac{\theta \ell^2-\ell^1}{1-\theta}, \ \delta_\theta \bar{\ell}:=|\delta_\theta\ell|+|\delta_\theta \widetilde{\ell}|,
\end{align*}
where $\ell=Y,y^{\tau},z^{\tau}$. Let $h$ be a positive constant only depending on the constants $\gamma_1,\gamma_2,\beta$ in Assumption \ref{ass quadratic'}. For any $M>0$, it is easy to check that 
\begin{equation}\label{esupY}\begin{split}
	&\E \left[\sup _{t \in[T-h, T]}\left|Y_t^1-Y_t^2\right|\right]=\E \left[\sup _{t \in[T-h, T]}\left|(1-\theta)(\delta_{\theta}Y_{t}+Y^1_t)\right|\right] \\
  \leq& (1-\theta)\E\left[\sup_{t \in[T-h, T]}\left|\delta_{\theta}Y_{t}\right|+\sup _{t \in[T-h, T]}\left|Y^1_t\right| \right]\leq (1-\theta)\E\left[\sup_{t \in[T-h, T]}\delta_\theta \bar{Y}_t+\sup _{t \in[T-h, T]}\left|Y^1_t\right| \right]\\
   \leq&(1-\theta)\left(\frac{1}{M} \sup _{\theta \in(0,1)} \E\left[\exp \left\{M \sup _{t \in[T-h, T]} \delta_\theta \bar{Y}_t\right\}\right]+\E\left[\sup _{t \in[T-h, T]}\left|Y_t^1\right|\right]\right).
  \end{split}
\end{equation}
As long as the following conditions holds
\begin{equation}\label{bddexpdetltay}
    \sup _{\theta \in(0,1)} \mathbf{E}\left[\exp \left\{M \sup _{t \in[T-h, T]} \delta_\theta \bar{Y}_t\right\}\right]<\infty,
\end{equation}
letting $\theta \to 1$ in \eqref{esupY}, we obtain that $Y^1=Y^2$ on $[T-h,T]$. Then, applying Itô's formula to $|Y^1-Y^2|^2$, we have $Z^1=Z^2,K^1=K^2$ on $[T-h,T]$. Repeating iteratively this procedure a finite number of times, we get the uniqueness on the whole interval. Therefore,  in the sequel, we prove that \eqref{bddexpdetltay} holds for appropriate $h$. 

\textbf{Step 1.} We show that, for any $t\in[0,T]$ and $p\geq 1$,
\begin{equation}\label{step 1}
    \begin{split}
        \exp\left(p\gamma\delta_\theta\bar{Y}_t\right)
    \leq \Bigg(\E_t\Bigg[\exp\bigg(p\gamma&\Big(\widetilde{\chi}+2(\gamma_1+\beta(T-t))\sup_{s\in[t,T]}|\delta_\theta Y_s|\\
    &+2(\gamma_2+\beta(T-t))\sup_{s\in[t,T]}\E\left[|\delta_\theta Y_s||\mathcal{G}_s\right]\Big)\bigg)\Bigg]\Bigg)^{2},
    \end{split}
\end{equation}
where 
\begin{align*}
\widetilde{\chi}&=2\chi+\sum_{i=1}^2\sup_{s\in[0,T]}|Y^i_s|, \ \ \zeta=max_{i=1,2} \sup_{t\in[0,T]}\E\left[\sup_{s\in[0,T]}|Y^i_s|\Big|\mathcal{G}_t\right],   \\
    \chi&=|\xi|+\alpha T+\sup_{s\in[0,T]}|h(s,0,0)|+2(\gamma_1+\beta T) \sum_{i=1}^2\sup_{s\in[0,T]}|Y^i_s|+2(\gamma_2+\beta T)\zeta.
\end{align*}

Simple calculation yields that $(\delta_\theta y^\tau,\delta_\theta z^\tau)$ is the solution to the following BSDE
\begin{equation}\label{eq18}
    \delta_\theta y^\tau=\delta_\theta \eta+\int_t^\tau \left(\delta_\theta f(s,\delta_\theta z^\tau_s)+\delta_\theta f_0(s)\right)ds-\int_t^\tau \delta_\theta z^\tau_s dB_s,
\end{equation}
where
\begin{align*}
    \delta_\theta \eta&=-\xi I_{\{\tau=T\}}+\frac{\theta (h\circ Y^1)_\tau-(h\circ Y^2)_\tau}{1-\theta}I_{\{\tau<T\}},\\
    \delta_\theta f_0(t)&=\frac{1}{1-\theta}\left((f\circ Y^1)(t,z^{2,\tau}_t)-(f\circ Y^2)(t,z^{2,\tau}_t)\right),\\
    \delta_\theta f(t,z)&=\frac{1}{1-\theta}\left(\theta(f\circ Y^1)(t,z^{1,\tau}_t)-(f\circ Y^1)(t,-(1-\theta)z+\theta z^{1,\tau}_t)\right).
\end{align*}
It follows from Assumption \ref{ass quadratic'} that 
\begin{equation}\begin{split}\label{eq18'}
    \delta_\theta \eta&\leq |\xi|+|h(\tau,0,0)|+\gamma_1(2|Y^1_\tau|+|\delta_\theta Y_\tau|)+\gamma_2\left(2\E\left[|Y^1_s|\big|\mathcal{G}_s\right]\big|_{s=\tau}+\E\left[|\delta_\theta Y_s|\big|\mathcal{G}_s\right]\big|_{s=\tau}\right),\\
    \delta_\theta f_0(t)&\leq \beta\left(|Y^1_t|+|\delta_\theta Y_t|+\E\left[|Y^1_t|\big|\mathcal{G}_t\right]+\E\left[|\delta_\theta Y_t|\big|\mathcal{G}_t\right]\right),\\
     \delta_\theta f(t,z)&\leq -(f\circ Y^1)(t,-z)\leq \alpha+\beta\left(|Y^1_t|+\E\left[|Y^1_t|\big|\mathcal{G}_t\right]\right)+\frac{1}{2}|z|^2. 
\end{split}\end{equation}
By Lemma \ref{technical lemma} and Estimate \eqref{eq18'}, for any $p\geq 1$, we have 
\begin{align*}
    &\exp\left(p\gamma (\delta_\theta y^\tau_t)^+\right)\\
    \leq &\E_t\left[\exp\left(p\gamma\left(\chi+(\gamma_1+\beta(T-t))\sup_{s\in[t,T]}|\delta_\theta Y_s|+(\gamma_2+\beta(T-t))\sup_{s\in[t,T]}\E\left[|\delta_\theta Y_s||\mathcal{G}_s\right]\right)\right)\right],
\end{align*}
which, together with \eqref{eq17} implies that 
\begin{equation}\label{eq19}
    \begin{split}
        &\exp\left(p\gamma (\delta_\theta Y_t)^+\right)\leq \esssup_{\tau\in\mathcal{T}_{t,T}}\exp\left(p\gamma (\delta_\theta y^\tau_t)^+\right)\\
    \leq &\E_t\left[\exp\left(p\gamma\left(\chi+(\gamma_1+\beta(T-t))\sup_{s\in[t,T]}|\delta_\theta Y_s|+(\gamma_2+\beta(T-t))\sup_{s\in[t,T]}\E\left[|\delta_\theta Y_s||\mathcal{G}_s\right]\right)\right)\right].
    \end{split}
\end{equation}
Similarly, we have
\begin{equation}\label{eq20}
    \begin{split}
        &\exp\left(p\gamma (\delta_\theta \widetilde{Y}_t)^+\right)\\
    \leq &\E_t\left[\exp\left(p\gamma\left(\chi+(\gamma_1+\beta(T-t))\sup_{s\in[t,T]}|\delta_\theta Y_s|+(\gamma_2+\beta(T-t))\sup_{s\in[t,T]}\E\left[|\delta_\theta Y_s||\mathcal{G}_s\right]\right)\right)\right].
    \end{split}
\end{equation}
Noting the fact that 
\begin{align*}
    (\delta_\theta Y)^-\leq (\delta_\theta \widetilde{Y})^++2|Y^2|, \ (\delta_\theta \widetilde{Y})^-\leq (\delta_\theta {Y})^++2|Y^1|,
\end{align*}
we obtain 
\begin{align*}
    &\exp\left(p\gamma\delta_\theta\bar{Y}_t\right)=\exp\left(p\gamma\left(|\delta_\theta {Y}_t|+|\delta_\theta \widetilde{Y}_t|\right)\right)\leq \exp\left(2p\gamma   \left((\delta_\theta Y_t)^++(\delta_\theta \widetilde{Y}_t)^++|Y^1_t|+|Y^2_t|\right)\right)\\
    \leq &\left(\E_t\left[\exp\left(p\gamma\left(\widetilde{\chi}+2(\gamma_1+\beta(T-t))\sup_{s\in[t,T]}|\delta_\theta Y_s|+2(\gamma_2+\beta(T-t))\sup_{s\in[t,T]}\E\left[|\delta_\theta Y_s||\mathcal{G}_s\right]\right)\right)\right]\right)^{2}.
\end{align*}

\textbf{Step 2.} We show that, for three constants $h\in(0,T]$, $p_1,p_2>1$ such that 
\begin{align*}
    4(\gamma_1+\gamma_2+2\beta h)<1, \ 4p_1(\gamma_1+\beta h)<1 \textrm{ and } 4p_1(\gamma_1+\beta h)<1,
\end{align*}
Estimate \eqref{bddexpdetltay} holds for $M=p\gamma$ with $p\geq 1$. 

In fact, recalling \eqref{step 1}, applying Doob's inequality and H\"{o}lder's inequality, for any $p\geq 1$, we have
\begin{align*}
    &\E\left[\exp\left(p\gamma \sup_{s\in[T-h,T]}\delta_\theta\bar{Y}_s\right)\right]\\
    \leq &4\E\left[\exp\left(2p\gamma\left(\widetilde{\chi}+2(\gamma_1+\beta h)\sup_{s\in[T-h,T]}|\delta_\theta Y_s|+2(\gamma_2+\beta h)\sup_{s\in[T-h,T]}\E\left[|\delta_\theta Y_s||\mathcal{G}_s\right]\right)\right)\right]\\
    \leq &
    4\left(\E\left[\exp\left(2p_3p\gamma\widetilde{\chi}\right)\right]\right)^{\frac{1}{p_3}}\left(\E\left[\exp\left(4p_1p\gamma(\gamma_1+\beta h)\sup_{s\in[T-h,T]}|\delta_\theta Y_s|\right)\right]\right)^{\frac{1}{p_1}}\\
    &\times \left(\E\left[\exp\left(4p_2p\gamma(\gamma_2+\beta h)\sup_{s\in[T-h,T]}\E\left[|\delta_\theta Y_s||\mathcal{G}_s\right]\right)\right]\right)^{\frac{1}{p_2}}\\
    \leq &
    4\left(\E\left[\exp\left(2p_3p\gamma\widetilde{\chi}\right)\right]\right)^{\frac{1}{p_3}}\left(\E\left[\exp\left(p\gamma\sup_{s\in[T-h,T]}|\delta_\theta Y_s|\right)\right]\right)^{4(\gamma_1+\beta h)}\\
    &\times \left(\frac{p_2}{p_2-1}\right) ^{p_2}\left(\E\left[\exp\left(4p_2p\gamma(\gamma_2+\beta h)\sup_{s\in[T-h,T]}|\delta_\theta Y_s|\right)\right]\right)^{\frac{1}{p_2}}\\
     \leq &
    4\left(\frac{p_2}{p_2-1}\right) ^{p_2}\left(\E\left[\exp\left(2p_3p\gamma\widetilde{\chi}\right)\right]\right)^{\frac{1}{p_3}}\left(\E\left[\exp\left(p\gamma\sup_{s\in[T-h,T]}|\delta_\theta Y_s|\right)\right]\right)^{4(\gamma_1+\gamma_2+2\beta h)},
\end{align*}
where $\frac{1}{p_1}+\frac{1}{p_2}+\frac{1}{p_3}=1$ and we have used Lemma \ref{technical lemma2} in the third inequality. Finally, we obtain that 
\begin{align*}
    \E\left[\exp\left(p\gamma \sup_{s\in[T-h,T]}\delta_\theta\bar{Y}_s\right)\right]\leq \left(4\left(\frac{p_2}{p_2-1}\right) ^{p_2}\left(\E\left[\exp\left(2p_3p\gamma\widetilde{\chi}\right)\right]\right)^{\frac{1}{p_3}}\right)^{\frac{1}{1-4(\gamma_1+\gamma_2+2\beta h)}}<\infty.
\end{align*}
The proof is complete. 
\end{proof}

\begin{remark}
    For the existence of the solution to conditional reflected BSDEs \eqref{conditional reflection} under Assumption \ref{ass quadratic'}, the usual construction is to define recursively a sequence of quadratic reflected BSDEs as follows: $Y^{(0)}_t=\E_t[\xi]$, $t\in [0,T]$ and for $m\geq 1$, $(Y^{(m)},Z^{(m)},K^{(m)})$ is the solution to reflected BSDE with terminal value $\xi$, generator $f^{(m)}$ and lower obstacle $h^{(m)}$, where $h^{(m)}_t=h(t,Y^{(m-1)}_t,\E[Y^{(m)-1}_t|\mathcal{G}_t])$ and $f^{(m)}$ can be chosen either of the following form
    \begin{align}
        &f^{(m)}(t,z)=f(t,Y^{(m-1)}_t,\E[Y^{(m-1)}_t|\mathcal{G}_t],z), \label{fm1} \\
        &f^{(m)}(t,y,z)=f(t,y,\E[Y^{(m-1)}_t|\mathcal{G}_t],z). \label{fm2}
    \end{align}
    Then, it remains to prove that $(Y^{(m)},Z^{(m)},K^{(m)})$ converges to $(Y,Z,K)$ under appropriate norm and $(Y,Z,K)$ is the solution to conditional reflected BSDE \eqref{conditional reflection}. The main problem here is the lack of the existence of solutions to the reflected BSDE with parameters $(\xi,f^{(m)},h^{(m)})$, where $\xi,h^{(m)}$ are clearly unbounded. To the best of our knowledge, the results in \cite{BY} are the most similar to the results we need here. It should be pointed out that the generator $g$ for the reflected BSDE in \cite{BY} should satisfy the following condition:
    \begin{align*}
        |g(t,y,z)|\leq \alpha+\beta|y|+\frac{\gamma}{2}|z|^2,
    \end{align*}
    where $\alpha,\beta,\gamma$ are positive constants.   Unfortunately, both $f^{(m)}$ defined in \eqref{fm1} and in \eqref{fm2} do not satisfy this condition.
\end{remark}

\section{Construction via penalization}

In this section, we present an alternative construction for the solutions to the conditional reflected BSDEs, namely the approximation via penalization, which is a frequently used method for classical reflected BSDEs. More precisely, the conditional reflected BSDE investigated in this section takes the following form
 \begin{align}\label{conditional reflection without Skorokhod}
       \begin{cases}
Y_t=\xi+\int_t^T f(s,Y_s,\E[Y_s|\mathcal{G}_s],Z_s)ds-\int_t^T Z_s dB_s+K_T-K_t, \\
Y_t\geq h(t,Y_t,\E[Y_t|\mathcal{G}_t]), \ t\in[0,T].
\end{cases}
   \end{align}
   Compared with the conditional reflected BSDEs in the previous section, here we do not impose the Skorokhod condition.  Since the Skorokhod condition usually ensures the minimality of the force $K$ aiming to push the solution upwards, an equivalent approach is to find the minimal solution $(Y,Z,K)\in \mathcal{S}^2\times\mathcal{H}^2\times\mathcal{S}^2_i$. That is, if $(Y',Z',K')\in\mathcal{S}^2\times\mathcal{H}^2\times\mathcal{S}^2_i$ is another solution to the conditional reflected BSDE \eqref{conditional reflection without Skorokhod}, we have $Y_t\leq Y'_t$, for any $t\in[0,T]$.

   It is worth pointing out that the construction via penalization relies heavily on the comparison theorem for conditional expectation BSDEs. Therefore, we need to propose the following additional structure of the generator $f$ and the loss function $h$.

\begin{assumption}\label{additional}
    \begin{itemize}
        \item[(i)] For each fixed $(t,y,z)$, the mapping $y'\rightarrow f(t,y,y',z)$ is non-decreasing.
        \item[(ii)] For any fixed $(t,y)$, the mapping $y'\rightarrow h(t,y,y')$ is non-decreasing.
        \item[(iii)] There exist $\eta\in L^2(\mathcal{F}_T)$ with $\eta\geq \xi$ and two measurable function $\Phi:[0,T]\times \Omega\times\mathbb{R}\times \mathbb{R}\times\mathbb{R}\rightarrow \mathbb{R}$ and $H:[0,T]\times \Omega\times\mathbb{R}\times \mathbb{R}\rightarrow \mathbb{R}$ such that for any $t,y,y',z$, we have
        \begin{align*}
            f(t,y,y',z)\leq \Phi(t,y,y',z), \ h(t,y,y')\leq H(t,y,y').
        \end{align*}
        Moreover, the conditional reflected BSDE \eqref{conditional reflection without Skorokhod} associated with terminal value $\eta$, generator $\Phi$ and loss function $H$ has a solution $(\bar{Y},\bar{Z},\bar{K})\in \mathcal{S}^2\times\mathcal{H}^2\times\mathcal{S}^2_i$.
    \end{itemize}
\end{assumption}

\begin{remark}\label{rem additional}
   (i) Suppose that the parameters $\eta,\Phi,H$ of the conditional reflected BSDE satisfy all the conditions as in Theorem \ref{main}. Clearly, the conditional reflected BSDE \eqref{conditional reflection without Skorokhod} associated with terminal value $\eta$, generator $\Phi$ and loss function $H$ has a solution. 

   \noindent (ii) Suppose that $f$ satisfies Assumption \ref{assGfh} (i) and for any $t,y,y'$, $h(t,y,y')\leq X_t$, where $X\in\mathcal{S}^2_c$. Then, similar with Theorem 5.2 in \cite{KKPPQ}, the conditional reflected BSDE \eqref{conditional reflection without Skorokhod} with terminal value $|\xi|+X_T$, generator $f$ and obsatcle $X$ has a solution $(\bar{Y},\bar{Z},\bar{K})\in\mathcal{S}^2_c\times \mathcal{H}^2\times\mathcal{S}^2_{ci}$ with 
   \begin{align*}
       \int_0^T (\bar{Y}_t-X_t)d\bar{K}_t=0.
   \end{align*}

   \noindent (iii) Suppose that $\eta,\Phi,H$ satisfy the following conditions:
   \begin{itemize}
       \item $\eta$ is a constant,
       \item $\Phi(t,y,y',z)=a^1_t y+a^2_t y'+a^3_t+g(t,z)$, where $a^1,a^2,a^3$ are deterministic and measurable functions such that $a^1,a^2$ are bounded and $a^3$ is square-integrable and $g$ satisfies Assumption \ref{assg} with $g(t,0)\equiv 0$,
       \item $H$ is a deterministic and continuous function such that $L(t,y)$ is strictly increasing and bi-Lipschitz in $y$, where $L(t,y)=y-H(t,y,y)$. 
   \end{itemize}
  We claim that the conditional reflected BSDE with parameters $\eta,\Phi,H$ has a solution. Actually, by Theorem 9 in \cite{BEH}, the following mean reflected BSDE
   \begin{align*}
       \begin{cases}
           Y_t=\eta+\int_t^T [(a^1_s+a^2_s)Y_s+a^3_s+g(s,Z_s)]ds-\int_t^T Z_s dB_s+K_T-K_t,\\
           \E[L(t,Y_t)]\geq 0, \ t\in[0,T] \textrm{ and } \int_0^T \E[L(t,Y_t)]dK_t=0
       \end{cases}
   \end{align*}
  has a unique solution $(Y,Z,K)\in \mathcal{S}^2_c\times \mathcal{H}^2\times \mathcal{S}^2_c$ with $K$ being a deterministic function.   It is easy to check $Z\equiv 0$ and thus, $Y$ is a deterministic function. Clearly, $(Y,0,K)$ is indeed the solution to the conditional reflected BSDE \eqref{conditional reflection without Skorokhod} with parameters $\eta,\Phi,H$.
\end{remark}

Suppose that Assumption \ref{assGfh} holds. For each $n\in \mathbb{N}$, by Theorem \ref{wellposedness for conditional BSDE},  the following conditional expectation BSDE parameterized by $n$ admits a unique pair of solution $(Y^n,Z^n)\in\mathcal{S}^2_c\times\mathcal{H}^2$:
\begin{equation}
    \begin{split}
        Y^n_t=\xi+\int_t^T f(s,Y^n_s,\E[Y^n_s|\mathcal{G}_s],Z^n_s)ds+n\int_t^T (Y^n_s-h(s,Y^n_s,\E[Y^n_s|\mathcal{G}_s]))^-ds-\int_t^T Z^n_s dB_s.
    \end{split}
\end{equation}
Now, we present the main result in this section. 

\begin{theorem}\label{construction via penalization}
    Let Assumptions \ref{assGfh} and \ref{additional} hold. Then, $(Y^n,Z^n,K^n)$ converges to $(Y,Z,K)\in\mathcal{S}^2\times\mathcal{H}^2\times\mathcal{S}^2_i$ in the sense that $K_t$ is the weak limit of $\{K^n_t\}_{n\in\mathbb{N}}$ in $L^2(\mathcal{F}_t)$ for each $t\in[0,T]$ and
   \begin{align*}
       &\lim_{n\rightarrow\infty}\E\left[\int_0^T|Y^n_s-Y_s|^2ds\right]=0,\\
       &\lim_{n\rightarrow\infty}\E\left[\int_0^T|Z^n_s-Z_s|^pds\right]=0, \textrm{ for any } p\in[1,2).
   \end{align*}
   Moreover, $(Y,Z,K)$ is the minimal solution to the  conditional reflected BSDE \eqref{conditional reflection without Skorokhod}.
\end{theorem}

\begin{remark}
    For the mean-field reflected BSDEs investigated in \cite{DEH}, the construction via penalization is also valid. The results (see Theorem 4.2 and Theorem 5.1 in \cite{DEH}) are stronger than ours, since the first component $Y$ and the compensating term $K$ are continuous and the Skorokhod condition is satisfied. However, it should be pointed out that the construction in \cite{DEH} is not merely a penalization scheme, but a recursive scheme of a penalization type (see Eqs. (4.3) and (5.3) in \cite{DEH}). Besides, the construction in \cite{DEH} requires that the loss function $h$ is non-decreasing in $y$, which is not needed in our setting. Moreover, if we can find the solution to conditional reflected BSDE associated with parameters $\eta,\Phi,H$ satisfying Assumption \ref{additional} (iii) (see Remark \ref{rem additional} for some typical cases), we could even drop the assumption that the Lipschitz constants $\gamma_1,\gamma_2$ of $h$ are sufficiently small, which is necessary both in \cite{DDZ,DEH,HMW} and in Theorem \ref{main}. 
\end{remark}

The proof of Theorem \ref{construction via penalization}  will be divided into several lemmas. We set 
\begin{align*}
    K^n_t=n\int_0^t (Y^n_s-h(s,Y^n_s,\E[Y^n_s|\mathcal{G}_s]))^-ds.
\end{align*}
For the approximation method applied in the classical reflected BSDEs, the first step is to provide the uniform estimate for the sequences $\{Y^n\}_{n\in\mathbb{N}}$, $\{Z^n\}_{n\in\mathbb{N}}$ and $\{K^n\}_{n\in\mathbb{N}}$ by using  It\^{o}'s  formula appropriately. However, this method is invalid in the present setting since the penalization term is far form being linear in the $y$-component. To overcome this problem, we will apply the comparison theorem for conditional expectation BSDEs (see Theorem \ref{comparison for conditional BSDE}). 
First, we show that the sequence $\{Y^n\}_{n\in\mathbb{N}}$ is non-decreasing in $n$, whose upper bound is $\bar{Y}$.

  \begin{lemma}\label{lem1}
      For any $n\in\mathbb{N}$, we have $Y^n\le Y^{n+1}\leq \bar{Y}$. 
  \end{lemma}  

  \begin{proof}
      Actually, for any $n\in\mathbb{N}$, $(Y^n,Z^n)$ is the solution to the conditional expectation BSDE with terminal value $\xi$ and generator $f^n$, where 
      \begin{align*}
          f^n(t,y,y',z)=f(t,y,y',z)+n(y-h(t,y,y'))^-.
      \end{align*}
      Clearly, $f^n$ is non-decreasing in $y'$ and $f^n(t,y,y',z)\leq f^{n+1}(t,y,y',z)$. It follows from Theorem \ref{comparison for conditional BSDE} that $Y^n\leq Y^{n+1}$. 

      Note that for any $n\in\mathbb{N}$, $(\bar{Y},\bar{Z})$ satisfies
      \begin{align*}
          \bar{Y}^n_t=\eta+\int_t^T \Phi^n(s,\bar{Y}_s,\E[\bar{Y}_s|\mathcal{G}_s],\bar{Z}_s)ds+\bar{K}_T-\bar{K}_t-\int_t^T\bar{Z}_sd B_s,
      \end{align*}
      where 
      \begin{align*}
          \Phi^n(t,y,y',z)=\Phi(t,y,y',z)+n(y-H(t,y,y'))^-.
      \end{align*}
      Applying Theorem \ref{comparison for conditional BSDE} again, we have $Y^n\leq \bar{Y}$.
  \end{proof}

  For any $t\in[0,T]$ and $n\in\mathbb{N}$, set
  \begin{align*}
      S^n_t=Y^n_t\wedge h(t,Y^n_t,\E[Y^n_t|\mathcal{G}_t]). 
  \end{align*}
  Clearly, for any $t\in[0,T]$, we have $Y^n_t\geq S^n_t$ and 
  \begin{align*}
      \int_t^T(Y^n_s-S^n_s)dK^n_s=n\int_t^T(Y^n_s-Y^n_s\wedge h(s,Y^n_s,\E[Y^n_s|\mathcal{G}_s]))(Y^n_s- h(s,Y^n_s,\E[Y^n_s|\mathcal{G}_s]))^-ds=0.
  \end{align*}
  Consequently, for any $t\in[0,T]$, we have
  \begin{equation}\label{YnKnSn}
      \int_t^T Y^n_sdK^n_s=\int_t^T S^n_s dK^n_s
  \end{equation}
 and $(Y^n,Z^n,K^n)$ can be viewed as the solution to the classical reflected BSDE with terminal value $\xi$, constant generator $\{f(s,Y^n_s,\E[Y^n_s|\mathcal{G}_s],Z^n_s)\}_{s\in[0,T]}$ and obstacle $S^n$. By Proposition 2.3 in \cite{KKPPQ}, we have 
 \begin{equation}\label{representation for Yn}
     Y^n_t=\esssup_{\tau\in\mathcal{T}_{t,T}}\E\left[\int_t^\tau f(s,Y^n_s,\E[Y^n_s|\mathcal{G}_s],Z^n_s)ds+S^n_\tau I_{\{\tau<T\}}+\xi I_{\{\tau=T\}}\Big|\mathcal{F}_t\right].
 \end{equation}
 The following lemma provides the uniform estimate for $\{Z^n\}_{n\in\mathbb{N}}$ and $\{K^n\}_{n\in\mathbb{N}}$, respectively.

  \begin{lemma}\label{lem2}
      There exists a constant $C$ independent of $n$, such that 
      \begin{equation}
          \E\left[\int_0^T|Z^n_s|^2ds\right]\leq C, \ \E\left[\int_0^T|f(s,Y^n_s,\E\left[Y^n_s|\mathcal{G}_s\right],Z^n_s)|^2ds\right]\leq C \textrm{ and } \E\left[|K^n_T|^2\right]\leq C.
      \end{equation}
  \end{lemma}

  \begin{proof}
      By Lemma \ref{lem1}, it is easy to check that 
      \begin{align}\label{bdd for Yn}
          \sup_{n\in\mathbb{N}}\E\left[\sup_{t\in[0,T]}|Y^n_t|^2\right]\leq \E\left[\sup_{t\in[0,T]}|Y^0_t|^2\right]+\E\left[\sup_{t\in[0,T]}|\bar{Y}_t|^2\right]<\infty.
      \end{align}
      Applying It\^{o}'s formula to $(Y^n_t)^2$ and noting \eqref{YnKnSn}, for any $\varepsilon>0$, we have 
      \begin{equation}\label{Ynito}\begin{split}
          &|Y^n_t|^2+\int_t^T|Z^n_s|^2 ds\\
          =&|\xi|^2+2\int_t^T Y^n_s f(s,Y^n_s,\E[Y^n_s|\mathcal{G}_s],Z^n_s)ds+2\int_t^T Y^n_s dK^n_s-2\int_t^T Y^n_sZ^n_s dB_s\\
          \leq &|\xi|^2+2\int_t^T |Y^n_s|\left(|f(s,0,0,0)|+L(|Y^n_s|+|\E[Y^n_s|\mathcal{G}_s]|+|Z^n_s|)\right)ds\\
          &+2K^n_T \sup_{s\in[0,T]}\{Y^n_s\wedge h(s,Y^n_s,\E[Y^n_s|\mathcal{G}_s])\}-2\int_t^T Y^n_sZ^n_s dB_s\\
          \leq &|\xi|^2+2\int_t^T |Y^n_s|\left(|f(s,0,0,0)|+L(|Y^n_s|+|\E[Y^n_s|\mathcal{G}_s]|)\right)ds+L^2\varepsilon^{-1}\int_t^T|Y^n_s|^2ds\\
          &+\varepsilon\int_t^T |Z^n_s|^2ds+\varepsilon|K^n_T|^2+\varepsilon^{-1} \left(\sup_{s\in[0,T]}\{Y^n_s\wedge h(s,Y^n_s,\E[Y^n_s|\mathcal{G}_s])\}\right)^2-2\int_t^T Y^n_sZ^n_s dB_s.
      \end{split}\end{equation}
      Set 
      \begin{align*}
          {h}^n_t=\gamma_1|Y^n_t|+\gamma_2\E[|Y^n_t||\mathcal{G}_t]+|h(t,0,0)|.
      \end{align*}
      Simple calculation yields that 
      \begin{align*}
          |Y^n_s\wedge h(s,Y^n_s,\E[Y^n_s|\mathcal{G}_s])-|Y^n_s|\wedge h^n_s|\leq |Y^n_s-|Y^n_s||+| h(s,Y^n_s,\E[Y^n_s|\mathcal{G}_s])-h^n_s|\leq 2(|Y^n_s|+h^n_s).
      \end{align*}
      Consequently,
      \begin{align*}
          |Y^n_s\wedge h(s,Y^n_s,\E[Y^n_s|\mathcal{G}_s])|\leq 3(|Y^n_s|+h^n_s). 
      \end{align*}
      Applying Doob's inequality and Eq. \eqref{bdd for Yn}, there exists a constant $C$ independent of $n$, such that 
      \begin{equation}\label{bdd for Sn}
          \E\left[\left(\sup_{s\in[0,T]}\{Y^n_s\wedge h(s,Y^n_s,\E[Y^n_s|\mathcal{G}_s])\}\right)^2\right]\leq C.
      \end{equation}
      Note that 
      \begin{align*}
          K^n_T=Y^n_0-\xi-\int_0^T f(s,Y^n_s,\E[Y^n_s|\mathcal{G}_s],Z^n_s)ds+\int_0^T Z^n_sdB_s.
      \end{align*}
      It is easy to check that 
      \begin{align}\label{bdd for Kn'}
          \E[|K^n_T|^2]\leq C+C\E\left[\int_0^T|Z^n_t|^2dt\right].
      \end{align}
      Taking expectations on both sides of \eqref{Ynito}, plugging in the inequality of Eqs. \eqref{bdd for Yn}, \eqref{bdd for Sn}, \eqref{bdd for Kn'} and choosing a sufficiently small $\varepsilon$, we obtain 
      \begin{align*}
          \E\left[\int_0^T|Z^n_t|^2dt\right]\leq C.
      \end{align*}
      Consequently, we have 
      \begin{align*}
          \E\left[\int_0^T|f(s,Y^n_s,\E\left[Y^n_s|\mathcal{G}_s\right],Z^n_s)|^2ds\right]\leq C \textrm{ and } \E\left[|K^n_T|^2\right]\leq C.
      \end{align*}
      The proof is complete.
  \end{proof}

  By Lemma \ref{lem1}, for any $t\in[0,T]$, we may define
  \begin{align*}
      Y_t:=\lim_{n\rightarrow\infty}Y^n_t.
  \end{align*}
  It is clear that 
  \begin{align}\label{bdd for Y}
          \E\left[\sup_{t\in[0,T]}|Y_t|^2\right]<\infty, \ \lim_{n\rightarrow\infty}\E\left[\int_0^T|Y^n_s-Y_s|^2ds\right]=0,
      \end{align} 
 which, together with Lemma \ref{lem1}, Lemma \ref{lem2} and Theorem 2.1 in \cite{P99}, implies the following result.

  \begin{lemma}\label{lem3}
      The process $Y$ is c\`{a}dl\`{a}g with representation
      \begin{align*}
          Y_t=\xi+\int_t^T g^0_sds+(K_T-K_t)-\int_t^T Z_sdB_s,
      \end{align*}
      where $g^0,Z\in\mathcal{H}^2$ are the weak limit of $\{f(s,Y^n_s,\E[Y^n_s|\mathcal{G}_s],Z^n_s),s\in[0,T]\}_{n\in\mathbb{N}}$, $\{Z^n\}_{n\in\mathbb{N}}$ in $\mathcal{H}^2$, respectively and $K\in S^2_i$ satisfies that $K_t$ is the weak limit of $\{K^n_t\}_{n\in\mathbb{N}}$ in $L^2(\mathcal{F}_t)$ for each $t\in[0,T]$. 
      Moreover,  for any $p\in[1,2)$, we have
      \begin{align}\label{converge of Z}
          \lim_{n\rightarrow\infty}\E\left[\int_0^T|Z^n_s-Z_s|^pds\right]=0.
      \end{align}
  \end{lemma}

  \begin{remark}
 It follows from the H\"{o}lder inequality and Eqs. \eqref{bdd for Y} and \eqref{converge of Z} that for any $p\in[1,2)$,
    \begin{align*}
        \lim_{n\rightarrow\infty}\E\left[\int_0^T |f(s,Y^n_s,\E[Y^n_s|\mathcal{G}_s],Z^n_s)-f(s,Y_s,\E[Y_s|\mathcal{G}_s],Z_s)|^pds\right]=0.
    \end{align*}
    By the uniqueness of the weak limit, $Y$ has the following decomposition
    \begin{align}\label{decomposition}
          Y_t=\xi+\int_t^T f(s,Y_s,\E[Y_s|\mathcal{G}_s],Z_s) ds+(K_T-K_t)-\int_t^T Z_sdB_s.
      \end{align}
  \end{remark}

  \begin{proof}[Proof of Theorem \ref{construction via penalization}]
      We firs show that for any $t\in[0,T]$, $Y_t\geq h(t,Y_t,\E[Y_t|\mathcal{G}_t])$. The inequality clearly holds for $t=T$. For $t<T$, it is a direct consequence of the following equation
     \begin{align*}
         \E\left[\int_t^T(Y_s-h(s,Y_s,\E[Y_s|\mathcal{G}_s]))^-ds\right]=\lim_{n\rightarrow\infty}\E\left[\int_t^T(Y^n_s-h(s,Y^n_s,\E[Y^n_s|\mathcal{G}_s]))^-ds\right]=0,
     \end{align*}
     where we have used Lemma \ref{lem2} in the last equality.

     It remains to prove that $(Y,Z,K)$ is the minimal solution. Let $(Y',Z',K')$ be the solution to conditional reflected BSDE \eqref{conditional reflection without Skorokhod}. For any $n\in\mathbb{N}$, it is clear that $(Y',Z',K')$ satisfies 
     \begin{align*}
         Y'_t=&\xi+\int_t^T f(s,Y'_s,\E[Y'_s|\mathcal{G}_s],Z'_s)ds-\int_t^T Z'_s dB_s+K'_T-K'_t\\
         &+n\int_t^T(Y'_s- h(s,Y'_s,\E[Y'_s|\mathcal{G}_s]))^-ds.
     \end{align*}
     By Theorem \ref{comparison for conditional BSDE}, we have $Y^n_t\leq Y'_t$, $t\in[0,T]$. Letting $n$ go to infinity, we obtain the desired result.
  \end{proof}

  The following lemma shows that a similar representation as \eqref{representation for Yn} also holds for the limit process $Y$ under some additional structure for $h$. 
  \begin{lemma}\label{lem4}
    Under Assumptions \ref{assGfh} and \ref{additional}, suppose that for any fixed $(t,y')$, the mapping $y\rightarrow h(t,y,y')$ is non-decreasing. For any $t\in[0,T]$, we have 
      \begin{equation}\label{representation for Y}
     Y_t=\esssup_{\tau\in\mathcal{T}_{t,T}}\E\left[\int_t^\tau f(s,Y_s,\E[Y_s|\mathcal{G}_s],Z_s)ds+h(\tau,Y_\tau,\E[Y_t|\mathcal{G}_t]|_{t=\tau}) I_{\{\tau<T\}}+\xi I_{\{\tau=T\}}\Big|\mathcal{F}_t\right].
 \end{equation}
  \end{lemma}

  \begin{proof}
      Recalling that $S^n_t=Y^n_t\wedge h(t,Y^n_t,\E[Y^n_t|\mathcal{G}_t])$, we define
      \begin{align*}
          G^n_t:=&\esssup_{\tau\in\mathcal{T}_{t,T}}\E\left[\int_t^\tau f(s,Y^n_s,\E[Y^n_s|\mathcal{G}_s],Z^n_s)ds+S^n_\tau I_{\{\tau<T\}}+\xi I_{\{\tau=T\}}\Big|\mathcal{F}_t\right]\\
          &-\esssup_{\tau\in\mathcal{T}_{t,T}}\E\left[\int_t^\tau f(s,Y_s,\E[Y_s|\mathcal{G}_s],Z_s)ds+S^n_\tau I_{\{\tau<T\}}+\xi I_{\{\tau=T\}}\Big|\mathcal{F}_t\right]
      \end{align*}
      and 
      \begin{align*}
          H^n_t:=\esssup_{\tau\in\mathcal{T}_{t,T}}\E\left[\int_t^\tau f(s,Y_s,\E[Y_s|\mathcal{G}_s],Z_s)ds+S^n_\tau I_{\{\tau<T\}}+\xi I_{\{\tau=T\}}\Big|\mathcal{F}_t\right].
      \end{align*}
      Clearly, we have $Y^n_t=G^n_t+H^n_t$. Note that 
      \begin{align*}
          |G^n_t|\leq \E\left[\int_0^T |f(s,Y^n_s,\E[Y^n_s|\mathcal{G}_s],Z^n_s)-f(s,Y_s,\E[Y_s|\mathcal{G}_s],Z_s)|ds\Big|\mathcal{F}_t\right].
      \end{align*}
     For any $p\in[1,2)$, applying Doob's inequality, H\"{o}lder's inequality, Fubini's theorem, Eqs. \eqref{bdd for Y} and \eqref{converge of Z} yields that 
     \begin{align*}
         \E\left[\sup_{t\in[0,T]}|G^n_t|^p\right]&\leq C\E\left[\int_0^T |f(s,Y^n_s,\E[Y^n_s|\mathcal{G}_s],Z^n_s)-f(s,Y_s,\E[Y_s|\mathcal{G}_s],Z_s)|^pds\right]\\
         &\leq C\left(\E\left[\int_0^T|Y^n_s-Y_s|^pds\right]+\E\left[\int_0^T|Z^n_s-Z_s|^pds\right]\right)\rightarrow 0, \textrm{ as } n\rightarrow\infty.
     \end{align*}
     Recall the fact that if for any $n\in\mathbb{N}$, $U^n$ and $U$ are c\`{a}dl\`{a}g and $U^n$ converges to $U$ increasingly, then the Snell envelope of $U^n$ converges to the Snell envelope of $U$ increasingly (see \cite{DHP} for more details). It is obvious that for any $t\in[0,T]$, $S^n_t\uparrow S_t:=Y_t\wedge h(t,Y_t,\E[Y_t|\mathcal{G}_t])$. It follows that $\{H^n\}_{n\in\mathbb{N}}$ converges increasingly to the following c\`{a}dl\`{a}g process $H$:
     \begin{align*}
         H_t:=\esssup_{\tau\in\mathcal{T}_{t,T}}\E\left[\int_t^\tau f(s,Y_s,\E[Y_s|\mathcal{G}_s],Z_s)ds+S_\tau I_{\{\tau<T\}}+\xi I_{\{\tau=T\}}\Big|\mathcal{F}_t\right].
     \end{align*}
     Recalling that for any $t\in[0,T]$, $Y_t\geq h(t,Y_t,\E[Y_t|\mathcal{G}_t])$, the proof is complete. 
  \end{proof}


 \section*{Acknowledgement}
    This work was supported  by the National Natural Science Foundation of China (No. 12301178), the Natural Science Foundation of Shandong Province for Excellent Young Scientists Fund Program (Overseas) (No. 2023HWYQ-049),  the Natural Science Foundation of Shandong Province (No. ZR2023ZD35),  the Fundamental Research Funds for the Central Universities and  the Qilu Young Scholars Program of Shandong University.

 \end{document}